\documentclass[a4paper,10pt]{article}
 \usepackage[utf8x]{inputenc}
 \usepackage[T1]{fontenc}
 \usepackage[normalem]{ulem}
 \usepackage[english]{babel}
 \usepackage[toc,page]{appendix} 
 \usepackage{amsfonts}
 \usepackage{graphicx}
  \usepackage{mathenv}
  \usepackage{amsmath}
  \usepackage[toc,page]{appendix}
 \usepackage{amsthm}
 \usepackage{pgf,tikz}
 \usetikzlibrary{arrows}
 \usetikzlibrary{shapes}
 \usetikzlibrary{plotmarks}
 \PrerenderUnicode{é}
  \author{ Pierre-Antoine Corre  \footnote{%
Laboratoire de Probabilités et Modèles Aléatoires, CNRS UMR 7599, Université Pierre et Marie Curie, Paris 6, Case courrier 188, 4 place Jussieu 75252 Paris Cedex 05, email: \texttt{pierre-antoine.corre@upmc.fr}}}
\newtheorem{theorem}{Theorem}[section]
\newtheorem*{thma}{Theorem (Chauvin, Rouault \cite{ChRoConnecting})}
\newtheorem*{thmb}{Theorem (Aïdékon \cite{aidekon2013convergence})}

\newtheorem*{thmd}{Theorem (Chauvin, Klein, Marckert and Rouault \cite{chauvin2005martingales})}
\newtheorem*{thme}{Theorem (Chauvin, and Drmota \cite{chauvin2006random})}
\newtheorem*{thmf}{Theorem (Drmota \cite{drmota2003analytic})}
\newtheorem*{thmg}{Theorem (Addario-Berry and Reed \cite{addario2009minima})}
 \newtheorem{lemma}[theorem]{Lemma}
  \newtheorem{proposition}[theorem]{Proposition}
  \newtheorem{corollary}[theorem]{Corollary}
   \newtheorem{Rem}{Remark}
    
  \title{Oscillations in the height of the Yule tree and application to the binary search tree}

\begin{document}
\definecolor{ffqqqq}{rgb}{1,0,0}
\definecolor{qqqqff}{rgb}{0,0,1}
\definecolor{xdxdff}{rgb}{0.49,0.49,1}
\maketitle
\begin{abstract}
For a particular case of a branching random walk with lattice support, namely the Yule branching random walk, we prove that the distribution of the centred maximum oscillates around a  distribution corresponding to a critical travelling wave in the following sense: there exist continuous functions $t \mapsto a_t$ and $x \mapsto \overline{\phi}(x)$ such that: 
 $$\lim_{t \rightarrow +\infty} \sup_{x \in \mathbb{R}} \vert \mathbb{P}(\overline{X}(t) \leq a_t +x )-\overline{\phi}(x- \{ a_t +x\})\vert=0,$$ where $\{x\}=x-\lfloor x \rfloor$ and $\overline{X}(t)$ is the height of the Yule tree. We also shows that similar oscillations occur for $\mathbb{E}\left(f(\overline{X}(t)-a_t)\right)$, when $f$ is in a large class of functions. This process is classically related to the binary search tree, thus yielding analogous results for the height and for the saturation level of the binary search tree. 
\end{abstract}

\section{Introduction} 
We denote by $\mathbb{N}$ the set $\lbrace 0, 1 , 2 ...\rbrace$ and by $\mathbb{N}^*$ the set $\mathbb{N} \setminus \lbrace 0 \rbrace$. 
The binary search tree $(\mathcal{T}_n)_{n \in \mathbb{N}^*}$ is a discrete time Markov process on the space of binary trees that can be constructed as follows: $\mathcal{T}_1$ is a tree made of a single leaf. Given $\mathcal T_n$, the next state $\mathcal{T}_{n +1}$ is obtained by uniformly choosing a leaf of $\mathcal{T}_n$ which we change into an internal node with two leaves attached. The binary search tree is a natural structure to store data and is related with the Quicksort algorithm. For a general reference on the binary search tree, see for instance the monograph of Mahmoud \cite{mahmoud1992evolution}. 

One can easily transform $(\mathcal{T}_n)_{n \in \mathbb{N}^*}$ into a continuous time Markov process simply by Poissonizing the jump times. More precisely, let $(N_t)_{t \in \mathbb{R}_+}$ be a pure birth process independent of $(\mathcal{T}_n)_{n \in \mathbb{N}^*}$ which jumps from state $n$ to state $n+1$ at rate $n$. 
Then 
\begin{equation}\label{equivyuletreecontbinarysearch}
(\mathcal{T}^c_t)_{t \in \mathbb{R}_+}:=(\mathcal{T}_{N_t})_{t \in \mathbb{R}_+}
\end{equation}
is called the Yule tree process. 
It is also a model of a random growing binary tree and is clearly a Markov process. Observe that we can give the following alternative description of the Yule tree process. At time $0$, the tree is reduced to a leaf. Each leaf lives for a random time with exponential distribution with parameter $1$ independent of the other leaves. When it dies, it is replaced by an internal node with two leaves. For more details and other constructions see e.g. \cite{chauvin2005martingales,reed2003height}.

Note that $N_t$ is the cardinal of $\mathcal{N}_t$, the set of leaves alive at time $t$, and that  if we introduce the stopping time 
\begin{equation}
\tau_n =\inf\lbrace t>0, N_t=n \rbrace, \label{tau}
\end{equation}
Equation \eqref{equivyuletreecontbinarysearch} yields:
\begin{equation}
\left(\mathcal{T}_n\right)_{n\in\mathbb{N}}=\left(\mathcal{T}^c_{\tau_n}\right)_{n\in\mathbb{N}} \label{embed1}
\end{equation}
so one can go from one model to the other.

From the Yule tree process we can make a branching random walk in the following way:
for a leaf $u \in \mathcal{N}_t$, let $X_u(t)$ be its height in the tree (or its generation). Define $(X(t))_{t\in \mathbb{R}_+}$, the measure valued process by:
\begin{equation}
X(t):=\sum_{u \in \mathcal{N}_t}\delta_{X_u(t)}, \; \forall t \in \mathbb{R}_+.
\end{equation}
$X$ is then simply a branching random walk in continuous time with lattice-integer support which we call the Yule branching random walk (sometimes it is called the Yule-time process as in \cite{ChRoConnecting}). Each particle lives for an exponential($1$) time and is then replaced by two daughter particles situated one unit of distance further.
Fix $0\leq t_0 \leq t_1 \leq t_2$. As usual, for $u \in \mathcal{N}_{t_1}$, we set $X_u(t_0)=X_v(t_0)$, where $v$ is the unique ancestor of $u$ in $\mathcal{N}_{t_0}$. Furthermore, for convenient we set $X_u(t_2)=X_u(t_1)$, even if $u$ is not alive at time $t_2$. We will also write $v<u$ when $v$ is the an ancestor of $u$. 
\\

The main focus of the present work is the study of the asymptotic behaviour, of the height (i.e. the highest generation of a leaf) and of the saturation level (i.e. the maximal level $l$ such that there are no leaves for all levels up to $l$) of the Yule tree and the binary search tree.  We point out that for the branching random walks we will talk about maximum and minimum rather than height and saturation level which are reserved to trees. 

 Our paper is organized as follows. In Section \ref{sectmainresult}, we state the main results. In Section \ref{sectprevresult}, we review the relevant literature and discuss our results. Section $\ref{SectYuletime}$ is dedicated to the proofs of results concerning the Yule branching random walk and Section \ref{sectbinasear} to the application of these results to the binary search tree. \\

\section{Main results \label{sectmainresult}}  
 Let us define $\overline{X}(t):=\max_{u \in \mathcal{N}_t} X_u(t)$ and $\underline{X}(t):=\min_{u \in \mathcal{N}_t} X_u(t).$ \\ 
 For $(x,t)\in \mathbb{R}\times \mathbb{R}^+$, set: 
\begin{equation}
\left\{
\begin{array}{l}
\overline{h}(x,t) := \mathbb{P}(\overline{X}(t) \leq x)= \mathbb{P}(\overline{X}(t) \leq \lfloor x \rfloor) \\
\underline{h}(x,t) := \mathbb{P}(\underline{X}(t) \geq x)=\mathbb{P}(\underline{X}(t) \geq \lceil x\rceil) \label{hmaxmin}
\end{array}.
\right.
\end{equation}
\begin{proposition}\label{propeqprin}
$\overline{h}$ and $\underline{h}$ solve the equation:
\begin{equation}
\partial_t h(x,t)= h^2(x-1,t)-h(x,t), (x,t)\in \mathbb{R}\times \mathbb{R}^+ .\label{eqnpri}
\end{equation}
\end{proposition}
The connection between the Yule branching random walk and Equation \eqref{eqnpri} is analogue to that between the branching Brownian motion and the F-KPP equation described in next section. We say that $x\mapsto \phi(x)$ is a travelling-wave solution of $\eqref{eqnpri}$ with speed $c\in \mathbb{R}$ if $(x,t) \mapsto \phi(x-ct)$ 
is a solution of \eqref{eqnpri}. 
One can easily check that $\phi$ is a travelling-wave solution of $\eqref{eqnpri}$ with speed $c\in \mathbb{R}$ if and only if it solves the differential equation:
\begin{equation}
c\phi^{\prime}(x)=\left(\phi(x)-\phi^2(x-1)\right), \quad x \in \mathbb{R}. \label{eqx}
\end{equation}

We will see in the next section that in the case of the branching Brownian motion and the F-KPP equation, the key result proved by Kolmogorov et al. in \cite{kolmogorov1937study} is that $t\mapsto \mathbb{P}(\overline{X}(t)\leq x+m_t)$ (where $m_t$ is an appropriate centring term) converges uniformly in $x$ to $\omega^*$, the critical traveling wave. 

In the case of the Yule branching random walk, the functions $\overline{h}$ and $\underline{h}$ defined in $\eqref{hmaxmin}$ are continuous in $t$ but piecewise constant in $x$, which implies that whatever the centring is, the distribution cannot converge. However, we see in the following theorem that the asymptotic distributions oscillate around the critical travelling-waves.  \\

For $\theta \in \mathbb{R}^*$, introduce 
\begin{equation}
c_{\theta}:=\frac{(2e^{\theta}-1)}{\theta}. \label{ctheta}
\end{equation}
We denote by $\theta^+$ (respectively by $\theta^-$) the largest (resp. the smallest) solution of $c_{\theta}=2e^{\theta}$. We also set $c^+:=c_{\theta^+}$ and $c^-:=c_{\theta^-}$.
Numerically, \cite{drmota2003analytic} we have: 
$$
\begin{cases}
c^+= 4.311...  \\
\theta^+=0.768...  \\
\end{cases} \qquad \qquad
 \begin{cases}
c^-= \; \;  \; 0.373...  \\
\theta^-=-1.678... \\
\end{cases}.     
$$

\begin{theorem}\label{thconvergenceYule}
Let $a_t=c^+ t-\frac{3 \log(t)}{2 \theta^+}$, $b_t=c^- t-\frac{3 \log(t)}{2 \theta^-}$ (recall that $\theta^-<0$) and $\{x\}=x-\lfloor x \rfloor$. \\
There exist a monotone travelling-wave solution $\overline{\phi}$ at speed $c^+$ and a monotone travelling-wave solution $\underline{\phi}$ at speed $c^-$ of  $\eqref{eqnpri}$  such that:
\begin{equation}
\lim_{t \rightarrow +\infty} \sup_{x \in \mathbb{R}} \vert \mathbb{P}(\overline{X}(t) \leq  a_t +x )-\overline{\phi}(x-\{a_t +x\})\vert=0 \label{converge}
\end{equation}
and:
\begin{equation}
\lim_{t \rightarrow +\infty} \sup_{x \in \mathbb{R}} \vert \mathbb{P}(\underline{X}(t) \geq  b_t +x )-\underline{\phi}(x -\{b_t +x\})\vert=0 . \label{converge2}
\end{equation}
\end{theorem}

Observe that \eqref{converge} is equivalent to:
\begin{equation}
\lim_{t \rightarrow +\infty} \sup_{x \in \mathbb{R}} \vert \mathbb{P}(\overline{X}(t) \leq  x )-\overline{\phi}(\lfloor x \rfloor-a_t)\vert=0
\end{equation}
and a similar formulation holds for \eqref{converge2}.
We will now extend Theorem $\ref{thconvergenceYule}$. For each piecewise continuous function $f$ we define $G_f$, $P_f$ and $r_f$ by:
\begin{equation}
\left\{
    \begin{array}{lll}
        G_f(t):=\mathbb{E}\left[f\left(\overline{X}(t)-a_t\right)\right], \quad t>0,  \\
        P_f(s):= \sum_{k \in \mathbb{Z}} f(k-s)\left( \overline{\phi}\left(k-s\right)-\overline{\phi}\left(k-1-s\right)\right), \quad s \in \mathbb{R},  \\ 
        r_f(t):= P_f(a_t), \quad t>0, \\
        \end{array}\label{gf1} 
\right. 
\end{equation}
when these functions are well-defined. Note that $P_f$ is clearly 1-periodic. Before stating the theorem, let us give some notations. We denote by $f^{(i)}$ the $i$th derivative of a function $f \in \mathcal{C}^{i}(\mathbb{R})$. For two real functions $f$ and $g$, we write:
$$f(x)=o_x\left(g(x)\right),$$
when there exists a function $\epsilon$ such that $\lim_{x \rightarrow + \infty} \epsilon(x)=0$ and such that $f(x)=\epsilon(x)g(x)$ and
$$f(x)=\underset{x \rightarrow \pm \infty}{o}\left(g(x)\right),$$
when the function $\epsilon$ also satisfies $\lim_{x \rightarrow - \infty} \epsilon(x)=0$.
\begin{theorem}\label{coroperio}
There exists $\delta>0$ such that for each piecewise continuous function $f$ which satisfies $f(x)=\underset{x \rightarrow \pm \infty}{o}(e^{\delta |x|})$ $G_f$, $P_f$ and $r_f$ are well-defined and we have:
\begin{equation}
\lim_{t \rightarrow +\infty} \lvert G_f(t)-r_f(t) \rvert=0. \label{firstpartcoroperio}
\end{equation}
Furthermore, if $f \in \mathcal{C}^k(\mathbb{R})$ satisfies $f^{(i)}(x)=\underset{x \rightarrow \pm \infty}{o}(e^{\delta |x|})$, for $0\leq i \leq k$, we have that $G_f,P_f,r_f \in \mathcal{C}^k(\mathbb{R^+})$ and:
\begin{equation}
\lim_{t \rightarrow +\infty} \lvert G_f^{(k)}(t)-r_f^{(k)}(t) \rvert =\lim_{t \rightarrow +\infty}\lvert G_f^{(k)}(t)-(c^+)^{k}P_f^{(k)}(a_t)\rvert=0. \label{secondpartcoroperio}
\end{equation}
\end{theorem}
Similar results for the minimum also hold. Observe that by taking, for $x$ fixed, $f$ defined by $f(y)=\mathbf{1}_{\lbrace y \leq  x \rbrace}$ in Theorem \ref{coroperio} we find the result of Theorem $\ref{thconvergenceYule}$ (if we omit the uniformity of the convergence).

It is interesting to see what Theorem \ref{coroperio} means in some particular cases. Set
\begin{equation}
\tilde{F}_t:=\mbox{Card}\lbrace u \in \mathcal{N}_t, \; X_u(t)=\overline{X}(t) \rbrace. \label{defFt}
\end{equation}
$\tilde{F}_t$ is the so-called "fringe" which has been studied by Roberts in \cite{roberts2010almost} and by Drmota in \cite{drmota2003analytic} (see next section for more details). Applying Theorem \ref{coroperio} with power functions  yields the following corollary.
\begin{corollary} \label{propuprim}
There exist three 1-periodic smooth functions $Q_1$,$Q_2$ and $Q_3$ such that:
 \begin{eqnarray}
&\mathbb{E}(\overline{X}(t)) &=a_t+Q_1(a_t)+o_t(1), \\
&\mathrm{Var}(\overline{X}(t)) &=Q_2(a_t)+o_t(1), \\
&\mathbb{E}(\tilde{F}_t) &=Q_3(a_t)+o_t(1).  \label{Ftperio}
\end{eqnarray}
\end{corollary}
To prove Equation \eqref{Ftperio}, we use the following lemma:
\begin{lemma}\label{lemmawt}
Let $w(t):=\mathbb{E}\left(\overline{X}(t)\right), \forall t\geq 0$. $w$ is a smooth function whose derivative is:
\begin{equation}
w^{\prime}(t)=\mathbb{E}\left(\tilde{F}_t\right), \forall t\geq 0.
\end{equation}
\end{lemma}



Let us now turn to the binary search tree. Using \eqref{equivyuletreecontbinarysearch} we show that Theorem $\ref{thconvergenceYule}$ translates into an analogous oscillation result for the asymptotic distributions of the height (i.e. the highest generation of a leaf) and the saturation level (i.e. the maximal level $l$ such that there are no leaves for all levels up to $l$) of a binary search tree. We call $\partial Z_{\infty}^+$ and $\partial Z_{\infty}^-$ the limit of the derivative martingales of the binary search tree defined in $\eqref{martderiveearbrebinaire}$, and for $K>0$:
\begin{equation}
\left\{
    \begin{array}{ll}
       \psi^+_K(x):=\mathbb{E}\left[\exp\left(-Ke^{-\theta^+ x} \partial Z_{\infty}^+\right)\right] \label{psiplus} \\
        \psi^-_K(x):=\mathbb{E}\left[\exp\left(-Ke^{-\theta^- x} \partial Z_{\infty}^-\right)\right] \\
       \end{array}
\right. .
\end{equation}
Let $H_n$ be the height  of a random binary search tree with $n$ nodes and $l_n$ its saturation level. 
\begin{theorem}\label{thconvergencearbre}
There exist $K^+,K^->0$ such that: 
\begin{equation}
\lim_{n \rightarrow +\infty} \sup_{x \in \mathbb{R}} \vert \mathbb{P}(H_n \leq \lfloor a_{\log(n)} +x \rfloor)-\psi^+_{K^+}(x-\{a_{\log(n)}+x\} )\vert=0, \label{convergearbre}
\end{equation}
and:
\begin{equation}
\lim_{n \rightarrow +\infty} \sup_{x \in \mathbb{R}} \vert \mathbb{P}(l_n \geq \lceil b_{\log(n)} +x \rceil)-\psi^-_{K^-}(x-\{b_{\log(n)}+x\} )\vert=0.
\end{equation}
\end{theorem}
We will further see in next section that Theorem \ref{thconvergencearbre} and results by Drmota allow us to obtain an analogue of Corollary \ref{propuprim}. Let $F_n$ be the number of particles at the highest position for a binary search tree with $n$ leaves. Note that $F_n$ is linked to $\tilde{F}_t$ by the relation $\tilde{F}_t=F_{N_t}$.
\begin{corollary} \label{corollarydrmotaperiodic}There exist three $1$-periodic functions $R_1$, $R_2$ and $R_3$ such that: 
 \begin{eqnarray}
& \mathbb{E}(H_n)&=a_{\log n}+R_1(a_{\log n})+o_n(1), \label{periodrmota}\\
& \operatorname{Var}(H_n)&=R_2(a_{\log n}) +o_n(1), \label{periodrmota2}\\
&\mathbb{E}(F_n)&=R_3(a_{\log n})+o_n(1).  \label{periodrmota3}
\end{eqnarray}
\end{corollary}

\section{Previous results and discussion\label{sectprevresult}} 
\subsection{Extremal particles in a branching process}\label{subsectiextrem}
 The position of the extremal particles in a branching process has been studied intensively. The case of the branching Brownian motion is the prototypical example. In that setting, it is then well known \cite{mckean1975application} that $p(x,t)=\mathbb{P}\left(\overline{X}(t) \leq x\right)$, where  $\overline{X}(t)$ is the position of the maximum at time $t$, solves the Fisher-Kolmogorov-Petrovskii-Piskunov (F-KPP) equation:
 \begin{equation}\label{Kpp}
\partial_t p(x,t)=\frac{1}{2}\partial_{x x}p(x,t)+\left(p^2(x,t)-p(x,t)\right).
\end{equation} 
Kolmogorov, Petrovskii and Piskunov \cite{kolmogorov1937study} show that, for a good centring term $m_t$,
\begin{equation}
 \lim_{t \rightarrow +\infty}p(x+m_t,t)=\omega^*(x), \label{convtravbbm}
 \end{equation}
where $\omega^* \in \mathcal{C}^2(\mathbb{R})$ is the unique function (up to a shift) such that $\tilde{p}(x,t):=\omega^*(x-\sqrt{2}t), \forall (x,t)\in \mathbb{R}\times \mathbb{R}_+$ is a solution of \eqref{Kpp}. 
One possible choice for $m_t$ is the median 
$m_t =\inf\{ x \in \mathbb{R} : p(t,x)=1/2\}$. Bramson \cite{bramson1983convergence} later shows famously that 
any valid centring term must be of the form
\begin{equation}
m_t =\sqrt 2 t -\frac3{2\sqrt 2} \log t +C +o(1).
\end{equation}
The study of the minimum (or equivalently of the maximum) of a branching random walk also has a long story. Let us mention some remarkable results on $M_n$, the minimum of a branching random walk. For simplicity, we consider a branching random walk satisfying suitable assumptions and which survives almost surely. Hammersley \cite{hammersley1974postulates} shows that $(M_n/n)$ converges almost surely to a constant $\gamma_1$. The almost sure behaviour of the minimum is refined by Hu and Shi \cite{hu2009minimal} who show that there exists a constant $\gamma_2>0$ such that: 
\begin{equation}
\liminf_{n \rightarrow + \infty} \frac{M_n-n\gamma_1}{\log n} =\frac{\gamma_2}{2}\mbox{ a.s. }, \quad \limsup_{n \rightarrow + \infty}\frac{M_n-n\gamma_1}{\log n}=\frac{3\gamma_2}{2}\mbox{ a.s. }.
\end{equation}
As far as the average of the minimum is concerned, Addario-Berry and Reed \cite{addario2009minima} show that:
\begin{equation}
\mathbb{E}(M_n)=\gamma_1 n+\frac{3\gamma_2}{2}\log n+O_n(1).
\end{equation}
They also show that the minimum centred around its mean is tight and that its distribution has exponential tails, that is there exist $C,\delta>0$ such that:
\begin{equation}\label{Addariointro}
\mathbb{P}(\lvert M_n-\mathbb{E}(M_n)\rvert\geq x)\leq Ce^{-\delta x}, \; \forall x \geq 0. 
\end{equation}
Aïdékon $\cite{aidekon2013convergence}$ proves a result similar to \eqref{convtravbbm} for a large class of branching random walks in the non-lattice case. He shows that there exists $K>0$ such that for all $x \in \mathbb{R}$:
\begin{equation}\label{Aidekonintro}
\lim_{n \rightarrow +\infty}\mathbb{P}\left(M_n\geq \gamma_1 n+\frac{3\gamma_2}{2}\log n+x\right)=\mathbb{E}\left(\exp\left(-Ke^{-x/\gamma_2}D_{\infty}\right)\right),
\end{equation}
where $D_{\infty}$ is the limit of the derivative martingale whose definition is given in \eqref{derivativeaidekon}. The results of Addario-Berry and Reed $\eqref{Addariointro}$ and of Aïdékon \eqref{Aidekonintro} will be useful for some proofs. They are therefore stated in more details in Appendix \ref{sectionresultsonmaxim}. 

Bramson, Ding and Zeitouni \cite{bramson2014convergence} give an alternative proof of Aïdékon's Theorem in a slightly less general case. In particular, they assume that the displacements of each offspring of a particle are independent. Furthermore, the authors claim that their method works in the lattice case. Note that, even if we consider a discretized version of the Yule branching random walk, since the displacements of the particles after a split are not independent in our framework, we cannot expect to use directly the method of \cite{bramson2014convergence}. 

 Let us now mention a work in the lattice case. Lifshits considers the following branching random walk in $\cite{lifshits2012cyclic}$. At time $n=0$, a particle is at $0$. At each time $n \in \mathbb{N}$, every particle produces two particles which are translated by $1$ from their parents with probability $0<p<1$ and by $-1$ with probability $1-p$. He proves that when $p>\frac{1}{2}$, the distribution of the centred maximum converges, but when $p=\frac{1}{2}$, oscillations as in $\eqref{converge}$ exist except that the centring term is not explicit (it involves the median). \\

 In that context, Theorem \ref{thconvergenceYule} provides an example of oscillations of the centred distribution of the maximum of a lattice branching random walk around a function (the critical travelling-wave) with an explicit centring term. We have already mentioned Aïdékon's result which shows that this kind of oscillations does not appears in the non-lattice case. It is interesting to see in Theorem \ref{coroperio} how this phenomena of oscillations extends to a large set of functions applied to $\overline{X}(t)-a_t$. Indeed, Theorem \ref{thconvergenceYule} only yields \eqref{firstpartcoroperio} for very restrictive classes of functions, for instance the class of continuous functions with compact support.
 
  Nevertheless, we point out that for a given function $f$, we cannot be sure, in general, that we have real oscillations in the sense that $P_f$ defined in \eqref{gf1} can be constant. In some specific cases, we can determinate whether $P_f$ is constant or not. 
 
 For instance, if $f$ is defined by $f(x)=\mathbf{1}_{x\leq y}$, then $P_f(x)=\overline{\phi}(y-\lbrace y-x \rbrace)$, and thus we have non-constant oscillations. Similarly, if we assume that $f$ is a non-constant 1-periodic function, then we have $P_f(x)=f(-x)$ and thus $P_f$ is also non-constant. However, for $f(x)=\overline{\phi}(x)+\overline{\phi}(x-1)$, we have that:
 $$P_f(x)=\sum_{k \in \mathbb{Z}} \overline{\phi}^2(k-x)-\overline{\phi}^2(k-x-1).$$
 Since $P_f$ is a telescoping sum, we have that $P_f(x)=1, \forall x \in \mathbb{R}$. 
 
 These examples are quite anecdotal but show that there is no trivial general answer to that issue. The examples of Corollary \ref{propuprim} seem more interesting. Whether $\mathbb{E}(\overline{X}(t)-a_t)$, Var($\overline{X}(t)$) and $\mathbb{E}(\tilde{F}_t)$ converge or not to a constant is still an open problem. The proof of Corollary \ref{propuprim} shows that $Q_3=c^+ + c^+Q^{\prime}_1$ which implies that $\mathbb{E}(\overline{X}(t)-a_t)$ converges to a constant if and only if $\mathbb{E}(\tilde{F}_t)$  also converges to a constant. Indeed, since $Q_1$ is $1$-periodic, $Q_1$ is constant if and only if $Q^{\prime}_1$ is constant. 
 
 Furthermore, the asymptotic behaviour of $\mathbb{E}(\tilde{F}_t)$ is especially interesting because it is related to a more general question, namely the convergence of the extremal point process. In the case of the branching Brownian motion, it has been shown independently by Aïdékon, Berestycki, Brunet and Shi  $\cite{aidekon2013branching}$ and Arguin, Bovier and Kistler $\cite{arguin2013extremal}$ that the extremal point process of branching Brownian motion converges and Madaule $\cite{madaule2011convergence}$ proved the analogous result for branching random walks with non-lattice support. The lattice case has not been dealt with and the behaviour of $\mathbb{E}(\tilde{F}_t)$ could shed a first light on this case. 
 
Analogous questions for the binary search tree will also be discussed in a next section.

\subsection{Yule generation process}

One of the key step in studying the Yule branching random walk is to switch our point of view by swapping the role of space and time. We thus introduce the Yule-generation process defined as follows. Let $\mathcal{M}_n$ be the set of the particles of the $n$th generation. For $n \in \mathbb{N}$ and for $u \in \mathcal{M}_n$ define 
\begin{equation}
T_u(n)=\inf\lbrace t>0, X_u(t)=n\rbrace \mbox{ and } T(n)=\sum_{u \in \mathcal{M}_n} \delta_{T_u(n)}. \label{nouveaubranch}
\end{equation} 
Observe that $T(n)$ is itself a branching random walk with discrete time and continuous spatial position with the following branching mechanism. 

At time $n=0$, a particle is at 0. At time $n=1$, the particle dies and gives birth to the point process: $\xi :=2 \delta_E$ where $E$ is distributed as an exponential variable with parameter $1$. We interpret $\xi$ as two particles, both in position $E$.  At each time $n$, the particles of the previous generation die and give birth to particles whose displacements from the parents are given by i.i.d. copies of $\xi$. 

This connection was established by Chauvin and Rouault $\cite{ChRoConnecting}$ and allowed them to prove the existence and uniqueness of travelling-waves. 

The proof of Theorem \ref{thconvergenceYule} relies on the correspondence between the Yule branching random walk and the Yule generation process defined in \eqref{nouveaubranch}. An advantage of such a change of point of view is that, as already mentioned, the non-lattice case is better understood. In particular, if we call $\underline{T}(n)=\inf \lbrace T_u(n),u \in \mathcal{M}_n \rbrace$,  Aïdékon's result can be transposed from the non-lattice case to our case thanks to the relation:
\begin{equation}
\mathbb{P}(\underline{T}(n) \leq t)=\mathbb{P}(\overline{X}(t) \geq n).\label{switch}
\end{equation}
Indeed, the event $\lbrace \underline{T}(n) \leq t \rbrace$ means that before $t$ there exists a particle whose generation is $n$, which is equivalent to the fact that the maximal generation of a particle at time $t$ is greater than $n$. To avoid confusions later on, we emphasize that if we take the complementaries of these events, \eqref{switch} yields:
\begin{equation}
\mathbb{P}(\underline{T}(n+1) \geq t)=\mathbb{P}(\overline{X}(t) \leq n).\label{switch2}
\end{equation}
An illustration of this correspondence is given below. The generation of the particle $u$ (the big circle) at time $t=4$ is $2$ and the reaching time of the second generation for $u$ is $T_u(2)$.

\begin{center}
\begin{tikzpicture}[line cap=round,line join=round,>=triangle 45, x=2cm,y=2cm]
\centering
\draw[->,color=black] (0,0) -- (4.5,0);
\foreach \x in {1,2,3,4}
\draw[shift={(\x,0)},color=black] (0pt,2pt) -- (0pt,-2pt) node[below] {\footnotesize $\x$};
\draw[->,color=black] (0,0) -- (0,3.5);
\foreach \y in {1,2,3}
\draw[shift={(0,\y)},color=black] (2pt,0pt) -- (-2pt,0pt) node[left] {\footnotesize $\y$};
\draw[color=black] (0pt,-10pt) node[right] {\footnotesize $0$};
\clip(-0.75,-0.5) rectangle (5,4);
\draw (0,0)-- (1.5,0);
\draw (1.5,1.05)-- (2.34,1.05);
\draw (3.68,1)-- (1.5,1);
\draw (2.34,2.05)-- (4.4,2.05);
\draw (3.32,1.95)-- (2.34,1.95);
\draw (3.68,2)-- (4.18,2);
\draw (3.68,2.1)-- (4.4,2.1);
\draw (3.32,3)-- (4.4,3);
\draw (4.4,3)-- (3.32,3);
\draw (3.32,3.05)-- (4.4,3.05);
\draw (4.18,3.1)-- (4.4,3.1);
\draw (4.18,3.15)-- (4.4,3.15);
\draw [dotted] (3.68,2.1)-- (3.68,0);
\draw (3.5,-0.15) node[anchor=north west] {$T_u(2)$};
\draw [dotted] (3.68,2)-- (0,2);
\draw (-0.75,2.3) node[anchor=north west] {$X_u(4)$};
\draw [dotted] (3.32,3.05)-- (3.32,0);
\draw (2.85,-0.15) node[anchor=north west] {$\underline{T}(3)$};
\draw [dotted] (3.32,3)-- (0,3);
\draw (-0.7,3.3) node[anchor=north west] {$\overline{X}(4)$};
\draw (-0.5,3.72) node[anchor=north west] {Generation};
\draw (4.50,0.15) node[anchor=north west] {Time};
\begin{scriptsize}
\draw [fill=xdxdff] (0,0) circle (1.4pt);
\draw [fill=qqqqff] (1.5,1) circle (1.4pt);
\draw [fill=qqqqff] (3.68,2) circle (2.3pt);
\draw [fill=qqqqff] (2.34,1.95) circle (1.4pt);
\draw [fill=qqqqff] (2.34,2.05) circle (1.4pt);
\draw [fill=qqqqff] (3.68,2.1) circle (1.4pt);
\draw [fill=qqqqff] (1.5,1.05) circle (1.4pt);
\draw [fill=qqqqff] (3.32,3) circle (1.4pt);
\draw [fill=qqqqff] (3.32,3.05) circle (1.4pt);
\draw [fill=qqqqff] (4.18,3.1) circle (1.4pt);
\draw [fill=qqqqff] (4.18,3.15) circle (1.4pt);
\draw [color=xdxdff] (3.68,0)-- ++(-1.5pt,-1.5pt) -- ++(3.0pt,3.0pt) ++(-3.0pt,0) -- ++(3.0pt,-3.0pt);
\draw [color=xdxdff] (0,2)-- ++(-1.5pt,-1.5pt) -- ++(3.0pt,3.0pt) ++(-3.0pt,0) -- ++(3.0pt,-3.0pt);
\draw [color=xdxdff] (3.32,0)-- ++(-1.5pt,-1.5pt) -- ++(3.0pt,3.0pt) ++(-3.0pt,0) -- ++(3.0pt,-3.0pt);
\draw [color=qqqqff] (0,3)-- ++(-1.5pt,-1.5pt) -- ++(3.0pt,3.0pt) ++(-3.0pt,0) -- ++(3.0pt,-3.0pt);
\end{scriptsize}
\end{tikzpicture}
\begin{figure}[h]
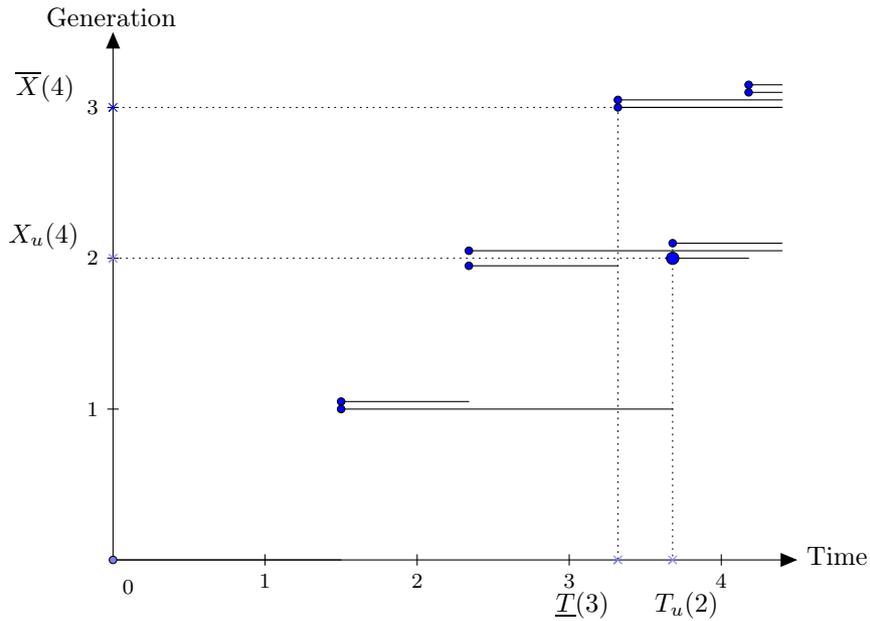

\caption{Correspondence between the Yule branching random walk and the Yule generation process}
\end{figure}
\end{center}

\subsection{Binary search tree}

The average and the variance of the height and of the saturation level have been studied by many authors. Initially, Robson $\cite{robson1979height,robson1982}$ proved that $\left(\mathbb{E}(H_n)/\log n\right)$ converges to a constant between $3.6$ and $c^+$ and in 1995, Devroye and Reed  $\cite{devroye1995variance}$ showed that $\mbox{Var}(H_n)=O_n((\log\log n)^2)$. These asymptotics have been improved $\cite{devroye1986note,devroye1987branching,drmota2001analytic,pittel1984growing}$ until Reed $\cite{reed2003height}$ and Drmota $\cite{drmota2003analytic}$ independently proved that 
\begin{equation}
\mathbb{E}(H_n)=a_{\log(n)}+O_n(1) \mbox{ and } \mbox{Var}(H_n)=O_n(1),
\end{equation}
where we recall that $a_t=c^+ t-\frac{3 \log(t)}{2 \theta^+}, \; \forall t \geq 0$.
 Moreover for $n \in \mathbb{N}$, let $\Upsilon_n$  be the generating function of $\left(\mathbb{P}(H_k\leq n)\right)_{k \in \mathbb{N}}$, that is:
 \begin{equation}
 \Upsilon_n(x):=\sum_{k=0}^{+\infty}\mathbb{P}(H_k\leq n)x^k, \; \forall x \in \mathbb{R}. \label{defyn}
 \end{equation}
 Drmota proves the convergence of the distribution of the height in the following sense.
\begin{thmf}
There exists a decreasing function $\Psi :\mathbb{R}^+ \rightarrow (0,1]$ with $\Psi(0)=1$ and 
$\lim_{x \rightarrow +\infty} \Psi(x)=0$ satisfying the integral equation
\begin{equation}
y\Psi(y/e^{(1/c^+)})=\int_0^y \Psi(z)\Psi(y-z)\mathrm{d}z
\end{equation}
such that
\begin{equation}
\lim_{k \rightarrow + \infty} \sup_{n \in \mathbb{N}} \left| \mathbb{P}(H_k \leq n)-\Psi\left(\frac{k}{\Upsilon_n(1)}\right)\right| =0. \label{equationconvdrmota}
\end{equation}
\end{thmf}
See also Drmota $\cite{drmota2009random}$ and Chauvin and Drmota $\cite{chauvin2006random}$. 
Moreover, Chauvin and Rouault \cite{ChRoConnecting} show that there exists a connection between $\Psi$, the function defined in Drmota's Theorem, and the derivative martingale of the binary search tree $\partial Z_{\infty}^+$  defined  in \eqref{martderiveearbrebinaire}. Indeed, there exists $K>0$ such that
\begin{equation}
\Psi(x)=\mathbb{E}\left[\exp\left(-Kx^{c^+-1}\partial Z_{\infty}^+\right)\right].\label{connectionPsi}
\end{equation} 
Observe that \eqref{connectionPsi} tells us that that $\psi^+_K(x)=\Psi(e^{-\frac{x}{c^+}})$, $\psi^+_K$ being defined in \eqref{psiplus}.

The main difference between Drmota's Theorem and Theorem \ref{thconvergencearbre} is that $\Upsilon_n(1)$ defined in \eqref{defyn} is implicit while $a_{\log n}$ is explicit. However, Drmota provides a good approximation of $\Upsilon_n(1)$ by showing that: \begin{equation}
\Upsilon_n(1)=e^{\frac{n}{c^+}+\frac{3\log n}{2(c^+-1)}+\kappa_n+o_n(1)},
\end{equation} 
where $(\kappa_n)$ is a bounded sequence such that $\kappa_{n+1}-\kappa_n \rightarrow 0$. He also proves that if there exists $\kappa_{\infty} \in \mathbb{R}$ such that the better asymptotic \begin{equation}
  \Upsilon_n(1)=e^{\frac{n}{c^+}+\frac{3\log n}{2(c^+-1)}+\kappa_{\infty}+o_n(1)}, \label{asymptyn}
  \end{equation}
holds then Corollary \ref{corollarydrmotaperiodic} also holds. By comparing Drmota's Theorem and Theorem \ref{thconvergenceYule} we can give an equivalent to $\Upsilon_n(1)$ and show that \eqref{asymptyn} holds thus yielding Corollary \ref{corollarydrmotaperiodic}. \\
  
 Drmota also gives more details about $F_n$. He shows in $\cite{drmota2004profile}$ that the oscillations of $(\mathbb{E}(F_n))$ around $c^+$ are at most of order $10^{-4}$. Moreover, in the same paper, he tells us that $(\mathbb{E}(F_n))$ is increasing until $n=100000$. If $(\mathbb{E}(F_n))$ were increasing on $\mathbb{N}$ we could easily prove that $(\mathbb{E}(F_n))$ converges. Whether or not $(\mathbb{E}(F_n))$ is increasing and whether or not  $(\mathbb{E}(F_n))$ converges are still open questions. 
 
 The issue of the almost sure behaviour of the height and of the saturation level for the binary search tree has been dealt with by Roberts in $\cite{roberts2010almost}$, where he shows, relying on \cite{hu2009minimal}, that: 
 \begin{equation}
 \frac{1}{2\theta^+}=\liminf_{n \rightarrow +\infty} \frac{c^+\log(n)-H_n}{\log\log(n)}<\limsup_{n \rightarrow +\infty} \frac{c^+\log(n)-H_n}{\log\log(n)}=\frac{3}{2\theta^+},
 \end{equation}
and the analogous result for the saturation level.  Roberts’ result yields a similar behaviour for the maximum and minimum of the Yule branching random walk. 

As far as the almost sure behaviour of $F_n$ is concerned, Roberts shows \cite{roberts2010almost} that:
\begin{equation}
\limsup_{n \rightarrow +\infty} F_n = +\infty.
\end{equation}

\section{Results on the Yule branching random walk}\label{SectYuletime}
\subsection{Proof of Proposition $\ref{propeqprin}$}
This is a simple transposition of McKean's proof for the branching Brownian motion \cite{mckean1975application}.
\begin{proof}[Proof of Proposition $\ref{propeqprin}$]
We recall that the first division time $\tau_2$ defined in \eqref{tau} has an exponential law with parameter 1. Consider $k \in \mathbb{Z}$ and $t\geq 0$. 
\begin{itemize}
\item If $k<0$ then $\overline{h}(k,t)=\mathbb{P}\left(\overline{X}(t)\leq k\right)=0, \forall t \geq 0$ and thus Equation $\eqref{eqnpri}$ is satisfied when $k<0$. 
\item If $k=0$, $\overline{h}(0,t)=\mathbb{P}\left(\overline{X}(t)\leq 0\right)=\mathbb{P}\left(\tau_2\geq t\right)=e^{-t}$.
Therefore
$$
\partial_t \overline{h}(0,t)=-e^{-t} =-\overline{h}(0,t) \nonumber =\overline{h}^2(-1,t) -\overline{h}(0,t).
$$
\item If $k>0$, we decompose the event $\lbrace \overline{X}(t)\leq k \rbrace$ in two parts depending on whether $\tau_2 \leq t$ or not. Since the event $\lbrace \tau_2 > t \rbrace$ is included in $\lbrace \overline{X}(t)\leq k \rbrace$, we have:
\begin{align*}
\mathbb{P}\left(\overline{X}(t)\leq k\right)&=\mathbb{P}\left(\tau_2 > t \right)+\mathbb{P}\left(\lbrace \overline{X}(t)\leq k \rbrace \cap \lbrace \tau_2 \leq t \rbrace\right) \\
&=e^{-t}+\mathbb{P}\left(\lbrace \overline{X}(t)\leq k \rbrace \cap \lbrace \tau_2 \leq t \rbrace\right).
\end{align*}
By strong Markov property, the event $\lbrace \overline{X}(t)\leq k \rbrace \cap \lbrace \tau_2 \leq t \rbrace$ is equal to $\lbrace \overline{X}^{(1)}(t-\tau_2)\leq k-1 \rbrace \cap\lbrace \overline{X}^{(2)}(t-\tau_2)\leq k-1 \rbrace \cap \lbrace \tau_2 \leq t \rbrace$, where $ \overline{X}^{(1)}$ and $\overline{X}^{(2)}$ are two independent copies of $\overline{X}$ which are also independent of $\tau_2$. Therefore,
\begin{align*}
\mathbb{P}\left(\overline{X}(t)\leq k\right)&=e^{-t}+\int_0^t e^{-s}\mathbb{P}\left(\overline{X}(t-s)\leq k-1\right)^2 \mathrm{d}s \\
&=e^{-t}+\int_0^t e^{-(t-s)}\mathbb{P}\left(\overline{X}(s)\leq k-1\right)^2 \mathrm{d}s \\
\Rightarrow e^t\mathbb{P}\left(\overline{X}(t)\leq k\right)&=1+\int_0^t e^{s}\mathbb{P}\left(\overline{X}(s)\leq k-1\right)^2 \mathrm{d}s.  
\end{align*}
By differentiating with respect to $t$, we obtain:
\begin{align*}
e^t(\partial_t \overline{h}(k,t)+\overline{h}(k,t))&=e^t \overline{h}^2(k-1,t) \\
\partial_t \overline{h}(k,t)&= \overline{h}^2(k-1,t)-\overline{h}(k,t),
\end{align*}
which is Equation \eqref{eqnpri}. 
\end{itemize}
We thus have proved that Equation \eqref{eqnpri} holds for $k \in \mathbb{Z}$. Since $\overline{h}(x,t)=\overline{h}(\lfloor x \rfloor,t)$, we have also proved that this equation holds for $x \in \mathbb{R}$.
\end{proof}
\subsection{Proof of Theorem $\ref{thconvergenceYule}$}
Since the proof of \eqref{converge} is analogous to that of \eqref{converge2}, we will just prove \eqref{converge}. We know that the critical travelling-waves are of the form $\phi_{K,\theta^+}$ (see \eqref{travcri}), where $K>0$. In order to simplify the notations, we will write $\phi_K$ instead of $\phi_{K,\theta^+}$ from now on. We thus want to prove that there exists $K>0$ such that:
\begin{equation}
\lim_{t \rightarrow +\infty} \sup_{x \in \mathbb{R}} \vert \mathbb{P}(\overline{X}(t) \leq\lfloor a_t +x \rfloor )-\phi_K(\lfloor a_t +x \rfloor -a_t)\vert=0 ,\label{re}
\end{equation}
where $a_t=c^+t-\frac{3\log(t)}{2 \theta^+}$. \\

Furthermore, remember that the branching random walk $(T(n))$ is obtained by switching time and space in the process $(X(t))$ (see $\eqref{nouveaubranch}$). The advantage of such a change of point of view is that $(T(n))$ is non-lattice and therefore, after some renormalizations, we can apply Aïdékon's Theorem (Appendix \ref{SubsecAid}) and finally use the relation $\eqref{switch}$ to prove $\eqref{re}$.

\begin{proof}[Proof of Theorem \ref{thconvergenceYule}]

Consider $T^{\prime}$ defined by:
\begin{equation} 
T^{\prime}(n)=\sum_{u \in \mathcal{M}_n} \delta_{(c^+-1)T_u(n)-\theta^+n}, \; \forall n \in \mathbb{N}.
\end{equation}
$T^{\prime}$ is a branching random walk which can also be described as follows. At time $n=0$, a particle is at $0$. For every $n\in \mathbb{N}$, each particle alive at time $n$ splits into two particles at time $n+1$. The displacement of the two new particles with respect to their parent is given by two independent random variables which have the same law as $(c^+-1)E-\theta^+$, where $E$ is an exponentially distributed random variable with parameter $1$.

 We will first ensure that $T^{\prime}$ satisfies Aïdékon's assumptions. Since the proofs of \eqref{assumpaidekon2} and \eqref{assumpaidekon3} are very close to that of Assumption \eqref{assumpaidekon1}, we will omit them. The first part of Assumption \eqref{assumpaidekon1} is obvious since there is two particles at each division. Now consider $E$ an exponentially distributed random variable with parameter $1$. We recall that the Laplace transform of $E$ is given by:
\begin{equation}
\mathbb{E}\left(e^{-\lambda E}\right)=\frac{1}{1+\lambda}. \label{laplacexp}
\end{equation}
The term of the second part of Assumption \eqref{assumpaidekon1} is:
\begin{align*}
\mathbb{E}\left(\sum_{u \in \mathcal{M}_1} e^{-T^{\prime}_u(1)}\right)&=2\mathbb{E}\left( e^{-\left[\left(c^+-1\right)E-\theta^+)\right]}\right) \\
&=2e^{\theta^+}\frac{1}{1+(c^+-1)}=2e^{\theta^+}\frac{1}{c^+}.
\end{align*}
By definition of $c^+$ and $\theta^+$, we have $c^+=2e^{\theta^+}$. Therefore:
\begin{equation}
\mathbb{E}\left(\sum_{u \in \mathcal{M}_1} e^{-T^{\prime}_u(1)}\right)=1
\end{equation}
and thus the second part of \eqref{assumpaidekon1} is proved.
By differentiating $\eqref{laplacexp}$ with respect to $\lambda$, we get:
\begin{equation}
\mathbb{E}\left(Ee^{-\lambda E}\right)=\frac{1}{(1+\lambda)^2}. \label{laplacexp2}
\end{equation}
This yields:
\begin{align*}
\mathbb{E}\left(\sum_{u \in \mathcal{M}_1} T^{\prime}_u(1)e^{-T^{\prime}_u(1)}\right)&=2\mathbb{E}\left(\left[\left(c^+-1\right)E-\theta^+)\right] e^{-\left[\left(c^+-1\right)E-\theta^+)\right]}\right) \\
&=2e^{\theta^+}\frac{c^+-1}{{(c^+)}^2}-\theta^+.
\end{align*}
Using the fact that $c^+=2e^{\theta^+}$ and $c^+=(2e^{\theta^+}-1)/\theta^+$, we get:
\begin{align*}
\mathbb{E}\left(\sum_{u \in \mathcal{M}_1} T^{\prime}_u(1)e^{-T^{\prime}_u(1)}\right)&=\frac{2e^{\theta^+}-1}{c^+}-\theta^+ \\
&=0.
\end{align*}
Hence, $T^{\prime}$ satisfies the assumptions of Aïdékon's Theorem. Now consider the derivative martingale $(D_n)$ defined in 
\eqref{derivativeaidekon} in the particular case of the branching random walk $(T^{\prime}(n))$:
\begin{align*}
D_n&=\sum_{u \in \mathcal{M}_n} T^{\prime}_u(n)e^{-T^{\prime}_u(n)} \\
&=\sum_{u \in \mathcal{M}_n} ((c^+-1)T_u(n)-\theta^+ n) e^{-[(c^+-1)T_u(n)-\theta^+ n]} \\
&=\theta^+ \partial W^{GEN}_n(\theta^+),
\end{align*}
where $\partial W^{GEN}_n(\theta^+)$ is defined in \eqref{martderiyulegen}. Furthermore, the relation \eqref{convderivgentotime} yields:
\begin{equation}
D_{\infty}=\theta^+ \partial W_{\infty}(\theta^+). \label{Dinfandderyulgen}
\end{equation}
For $n \in \mathbb{N}$ and $x \in \mathbb{R}$, we define: 
\begin{equation}
t_n=\frac{n}{c^+}+\frac{3\log(n)}{2 c^+\theta^+} \quad \mbox{ and } \quad d_n(x)=\mathbb{P}\left(\underline{T}(n)\geq t_n-\frac{x}{c^+}\right). \label{deftn}
\end{equation}
After straightforward computations, we can show that
$$
d_n(x)=\mathbb{P}\left(\underline{T}^{\prime}(n)\geq \frac{3}{2}\log n+x\right),
$$
where $\underline{T}^{\prime}(n)=\min \lbrace T^{\prime}_u(n), u \in \mathcal{M}_n \rbrace$. Aïdékon's result \eqref{eqthAidekon} and Equation \eqref{Dinfandderyulgen} yield that there exists $K>0$ such that:
\begin{equation}
\lim_{n \rightarrow +\infty} d_n(x) = \phi_K(x), \forall x \in \mathbb{R}. \label{convunf}
\end{equation}
Moreover, for each $n \in \mathbb{N}$, the fact that $d_n$ is a cumulative distribution function and the form of $\phi_K$ \eqref{travcri} ensure that $(d_n)$ and $\phi_K$ satisfy the assumptions of Theorem \ref{polyatheorem}. Therefore:
\begin{equation}
\lim_{n \rightarrow +\infty} \sup_{x \in \mathbb{R}} \left| d_n(x) - \phi_K(x)\right|=0. \label{convunf2}
\end{equation}
We will now use $\eqref{switch2}$ to prove $\eqref{re}$. We start by proving that there exists $K>0$ such that for all $x \in \mathbb{R}$:
\begin{equation}
\lim_{t \rightarrow +\infty} \vert \mathbb{P}(\overline{X}(t) \leq\lfloor a_t + x\rfloor )-\phi_{K}(\lfloor a_t + x\rfloor -a_t)\vert=0. \label{moins}
\end{equation}
We first rewrite $\mathbb{P}(\overline{X}(t) \leq \lfloor a_t + x\rfloor)$ by using $\eqref{switch2}$:
\begin{align}
\mathbb{P}(\overline{X}\left(t\right) \leq \lfloor a_t + x\rfloor) &= \mathbb{P}\left(\underline{T}(\lfloor a_t + x\rfloor +1) \geq t\right) \nonumber \\
&=\mathbb{P}\left(\underline{T}(\lfloor a_t + x\rfloor+1) \geq t_{\lfloor a_t + x\rfloor+1}+t-t_{\lfloor a_t + x\rfloor+1}\right) \nonumber \\
&=\mathbb{P}\left(\underline{T}(\lfloor a_t + x\rfloor+1) \geq t_{\lfloor a_t + x\rfloor+1}-\frac{c^+( t_{\lfloor a_t + x\rfloor+1}-t)}{c^+} \right) \nonumber \\
&=d_{\lfloor a_t + x\rfloor+1}\left(c^+\left( t_{\lfloor a_t + x\rfloor+1}-t\right)\right). \label{argg}
\end{align}
A straightforward computation gives:
$$c^+( t_{\lfloor a_t + x\rfloor+1}-t)=\lfloor a_t + x\rfloor-a_t+1+\frac{3\log(c^+)}{2\theta^+}+o_t(1),$$
where $o_t(1)=\frac{3}{2c^+\theta^+}\log\left(\frac{\lfloor a_t +x\rfloor+1}{a_t+x+1}\right)$. By the triangle inequality:
\begin{align*}
&\left\lvert d_{\lfloor a_t + x\rfloor+1}(c^+( t_{\lfloor a_t + x\rfloor+1}-t))-\phi_K \left(\lfloor a_t + x\rfloor-a_t+1+\frac{3\log(c^+)}{2\theta^+}\right) \right\rvert \\
&\; \leq \quad  \left\lvert d_{\lfloor a_t + x\rfloor+1}(c^+( t_{\lfloor a_t + x\rfloor+1}-t))-\phi_K \left(c^+( t_{\lfloor a_t + x\rfloor+1}-t)\right) \right\rvert \\
&  \quad  \quad \; + \left\lvert \phi_K \left(c^+( t_{\lfloor a_t + x\rfloor+1}-t)\right) -\phi_K \left(\lfloor a_t + x\rfloor-a_t+1+\frac{3\log(c^+)}{2\theta^+}\right) \right\rvert,
\end{align*}
where $K$ is fixed by $\eqref{convunf}$. The first term on the right hand side goes to $0$ by uniform convergence of $(d_n)$ to $\phi_K$, see $\eqref{convunf2}$, and the second goes to $0$ by uniform continuity of $\phi_K$ (by \eqref{travcri}  it is a continuous function with  finite limits as $x \rightarrow \pm \infty$). 

We thus have proved that for all $x \in \mathbb{R}$:
$$ \lim_{t \rightarrow + \infty} \left\lvert \mathbb{P}(\overline{X}(t) \leq \lfloor a_t + x\rfloor) -\phi_K \left(\lfloor a_t + x\rfloor-a_t+1+\frac{3\log(c^+)}{2\theta^+}\right) \right\rvert=0,$$
which is $\eqref{moins}$ (since the travelling-wave is defined up to an additive translation in the argument, we can incorporate $1+\frac{3\log(c^+)}{2\theta^+}$ into $K$). 

With the help of $\eqref{moins}$, we will now prove \eqref{re}. For that purpose, we start by proving that there exist $t_0,x_0,x_1 \in \mathbb{R}$, $x_0<x_1$ such that: 
\begin{equation}
\sup_{x \in ]-\infty,x_0] \cup [x_1,+\infty[} \vert \mathbb{P}(\overline{X}(t) \leq\lfloor a_t + x\rfloor)-\phi_{K}(\lfloor a_t + x\rfloor-a_t)\vert<\epsilon, \; \forall t>t_0. \label{supext}
 \end{equation}
 First observe that:
 $$x \leq\lfloor a_t+x \rfloor - a_t<x+1, \; \forall x \in \mathbb{R}, \; \forall t>0.$$ 
 Furthermore, $\phi_K$ is increasing and tends to $0$ when $x$ goes to $-\infty$. Therefore, for $\epsilon>0$, there exists $x_0\in \mathbb{R}$ such that: 
\begin{equation}
0<\phi_{K}( \lfloor a_t+x \rfloor - a_t)<\epsilon/3, \forall t>0, \forall x\leq x_0.\label{epsilon3part1aidekonequivproof}
  \end{equation}
  Equation  $\eqref{moins}$ yields the existence of $t_1>0$ such that:
\begin{equation}
\vert \mathbb{P}(\overline{X}(t) \leq\lfloor a_t +x_0\rfloor)-\phi_{K}(\lfloor a_t +x_0\rfloor-a_t)\vert<\epsilon/3, \; \forall t>t_1. \label{epsilon3part2aidekonequivproof}
\end{equation}
 The functions $x \mapsto \mathbb{P}(\overline{X}(t) \leq \lfloor a_t+x \rfloor)$ is non-decreasing and tends to $0$ when $x$ goes to $-\infty$. Therefore, Equations \eqref{epsilon3part1aidekonequivproof} and \eqref{epsilon3part2aidekonequivproof} yield
\begin{equation} 
\vert \mathbb{P}(\overline{X}(t) \leq\lfloor a_t +x\rfloor)\vert<2\epsilon/3, \forall x\leq x_0, \forall t>t_1.\label{epsilon3part3aidekonequivproof}
 \end{equation}
Moreover, using the triangle inequality and with the help of \eqref{epsilon3part1aidekonequivproof} and \eqref{epsilon3part3aidekonequivproof}, we get:
\begin{equation} 
\vert \mathbb{P}(\overline{X}(t) \leq\lfloor a_t +x\rfloor)-\phi_{K}(\lfloor a_t +x\rfloor-a_t)\vert<\epsilon, \forall x\leq x_0, \forall t>t_1.
 \end{equation}
We can prove in a similar way that there exists $t_2>0$ and $x_1>x_0$ such that:
\begin{equation} 
\vert \mathbb{P}(\overline{X}(t) \leq\lfloor a_t +x\rfloor)-\phi_{K}(\lfloor a_t +x\rfloor-a_t)\vert<\epsilon, \forall x\geq x_1, \forall t>t_2.
 \end{equation}
By taking $t_0=\max \lbrace t_1, t_2 \rbrace$, we thus have \eqref{supext}.
Moreover, since $[x_0,x_1] \cap \mathbb{Z}$ is finite, we have that:
\begin{equation}
\lim_{t \rightarrow +\infty} \sup_{x \in [x_0,x_1] \cap \mathbb{Z}} \vert \mathbb{P}(\overline{X}(t) \leq\lfloor a_t + x\rfloor)-\phi_{K}(\lfloor a_t + x\rfloor-a_t)\vert=0. \label{supintZ}
\end{equation}
For  fixed $t$, $x \mapsto \vert \mathbb{P}(\overline{X}(t) \leq\lfloor a_t + x\rfloor)-\phi_{K}(\lfloor a_t + x\rfloor-a_t)\vert$ is constant on each interval of the form: $[\lfloor a_t \rfloor-a_t  +k, \lfloor a_t  \rfloor -a_t+k+1[$, and then for fixed $t$ 
\begin{align}
&\sup_{x \in [x_0,x_1] \cap \mathbb{Z}} \vert \mathbb{P}(\overline{X}(t) \leq\lfloor a_t + x\rfloor)-\phi_{K}(\lfloor a_t + x\rfloor-a_t)\vert \nonumber \\
&\qquad =\sup_{x \in [x_0,x_1]} \vert \mathbb{P}(\overline{X}(t) \leq\lfloor a_t + x\rfloor)-\phi_{K}(\lfloor a_t + x\rfloor-a_t)\vert. \label{coresZR}
\end{align}
Combining $\eqref{supext}$, $\eqref{supintZ}$ and $\eqref{coresZR}$ we obtain:
\begin{equation}
\lim_{t \rightarrow +\infty} \sup_{x \in \mathbb{R}} \vert \mathbb{P}(\overline{X}(t) \leq\lfloor a_t +x \rfloor )-\phi_{K}(\lfloor a_t +x \rfloor -a_t)\vert=0 .
\end{equation}
\end{proof}

\subsection{Proof of Theorem $\ref{coroperio}$}
To prove Theorem $\ref{coroperio}$,  we first need to establish estimates on the tails of the distribution of the maximum to determine a class of functions $f$ sufficiently large for which $G_f$ and $P_f$ $\eqref{gf1}$ are well defined. Again, we will work with the Yule generation process. More precisely, we rely on a result due to Addario-Berry and Reed $\cite{addario2009minima}$ (see Appendix \ref{Subsecaddario}) which shows that the maximum of a branching random walk with i.i.d. displacements has exponential tails. 

\begin{proof}[Proof of Theorem \ref{coroperio}]
Although the process $T$ does not have independent displacements, it is very easy to modify it to fit the hypothesis of Addario-Berry and Reed's Theorem (see Appendix \ref{Subsecaddario}). We can consider ($T^{(2)}(n)$) defined by: 
\begin{equation}
T^{(2)}(n)=\frac{1}{2}\sum_{u \in \mathcal{M}_{n+1}} \delta_{T_u(n+1)-\tau_2}, \quad \forall n \geq 0, \label{Tprim}
\end{equation}
where $\tau_2$ is defined in $\eqref{tau}$. Note that $(T^{(2)}(n))$ is just a translation of $(T(n))$ where we remove one element of each couple of particles. The branching random walk $(T^{(2)}(n))$ can also be described as in Appendix \ref{Subsecaddario} by taking $B=2$ and $Y=E$, where $E$ is an exponential random variable with parameter $1$. Furthermore, $(T^{(2)}(n))$ satisfies the assumptions of Addario-Berry and Reed's Theorem. Indeed, Assumptions 1,2 and 3 are clearly satisfied. By taking $\lambda_1=-c^+\theta^+$, we can also show, after straightforward computations, that Assumption 4 is satisfied.  

 Let $\underline{T}^{(2)}=\min \lbrace T^{(2)}_u(n), u \in \mathcal{M}_n \rbrace $ and recall that $t_n$ is defined by $\eqref{deftn}$. Since the branching random walk $T^{(2)}$ almost surely survives, Addario-Berry and Reed's Theorem yields that there exists a bounded sequence $\rho$ such that:
\begin{equation}
\mathbb{E}\left(\underline{T}^{(2)}(n)\right)=t_n+\rho_n, \; \forall n \geq 1
\end{equation}
and that there exist $C,\alpha>0$  such that for all $n \geq 1$ and $x \in \mathbb{R}^+$:
\begin{equation}
\mathbb{P}\left(|\underline{T}^{(2)}(n)-t_n-\rho_n|\geq x\right) \leq Ce^{-\alpha x}. \label{tailTminmintnprovisoire}
\end{equation}
Since for each $x \in \mathbb{R}^+$, $\alpha \mapsto e^{-\alpha x}$ is non-increasing we can suppose that $\alpha<1$.
Furthermore, the relation $\eqref{Tprim}$ clearly implies:
\begin{equation}
\underline{T}(n) = \underline{T}^{(2)}(n-1)+\tau_2 \label{encadreTparTprime}
\end{equation}
and that $\tau_2$ is independent of $T^{(2)}$. Therefore, for $n\geq 2$ and $x\in \mathbb{R}$ we get by the triangular inequality:
\begin{align*}
&\mathbb{P}\left(\left.\left|\underline{T}^{(2)}(n-1)-t_{n-1}-\rho_{n-1}\right|\geq x-\tau_2 \right|  \tau_2\right) \leq Ce^{-\alpha (x-\tau_2)} \\
\Rightarrow \quad &\mathbb{P}\left(\left.\left|\underline{T}^{(2)}(n-1)-t_{n-1}-\rho_{n-1}\right|+\tau_2\geq x\right|\tau_2\right) \leq Ce^{\alpha \tau_2}e^{-\alpha x} \\
\Rightarrow \quad &\mathbb{P}\left(\left.\left|\underline{T}^{(2)}(n-1)+\tau_2-t_{n-1}-\rho_{n-1}\right|\geq x\right|\tau_2\right) \leq Ce^{\alpha \tau_2}e^{-\alpha x} \\
\Rightarrow  \quad &\mathbb{P}\left(\left.\left|\underline{T}(n)-t_n+t_n-t_{n-1}-\rho_{n-1}\right|\geq x\right|\tau_2\right) \quad \; \leq Ce^{\alpha \tau_2}e^{-\alpha x}. \\
\end{align*}
The sequence $(|t_n-t_{n-1}-\rho_{n-1}|)$ is bounded, therefore there exists $M>0$ such that:
\begin{align*}
&\mathbb{P}\left(|\underline{T}(n)-t_n|-M\geq x |\tau_2\right) \leq Ce^{\alpha \tau_2}e^{-\alpha x} \\
\Rightarrow \quad &\mathbb{P}\left(|\underline{T}(n)-t_n|\geq x|\tau_2\right) \quad \leq Ce^{\alpha M}e^{\alpha \tau_2}e^{-\alpha x}. \\
\end{align*}
We have assumed that $\alpha<1$ and thus $\mathbb{E}\left(e^{\alpha \tau_2}\right)<+\infty$. This yields that there exists $C_2>0$ such that:
\begin{equation}
\mathbb{P}\left(|\underline{T}(n)-t_n|\geq x\right) \leq C_2e^{-\alpha x}. \label{tailTminmintn}
\end{equation}
Fix $\eta>0$. Since $a$ is not bounded in the neighbourhood of $0$, we will now choose our times $t$ in the set $[\eta,+\infty)$. With the help of $\eqref{switch}$ and \eqref{tailTminmintn}, we now want to obtain that there exist $C_3>0$ and $\alpha^{\prime}>0$ such that for all $t \geq \eta$ and all $k \in \mathbb{N}$
\begin{equation}
\mathbb{P}\left(|\overline{X}(t)-\lfloor a_t \rfloor |\geq k\right) \leq C_3e^{-\alpha^{\prime} k}. \label{asymptailH}
\end{equation}
For that purpose, it suffices to prove that there exist $C_4,C_5>0$ such that:
\begin{equation}
\mathbb{P}\left(\overline{X}(t)-\lfloor a_t \rfloor \geq k\right) \leq C_4e^{-\alpha^{\prime} k} \mbox{ and } \mathbb{P}\left(\overline{X}(t)-\lfloor a_t \rfloor \leq -k\right) \leq C_5e^{-\alpha^{\prime} k}. \label{distingxatk}
\end{equation}
Let us prove the first part of \eqref{distingxatk}. By taking $n=\lfloor a_t \rfloor + k$ in $\eqref{tailTminmintn}$, we get:
\begin{align*}
&\mathbb{P}\left(|\underline{T}(\lfloor a_t \rfloor + k)-t_{\lfloor a_t \rfloor + k}|\geq x\right) \leq C_2e^{-\alpha x} \\
\Rightarrow \quad &\mathbb{P}\left(\underline{T}(\lfloor a_t \rfloor + k)\leq t_{\lfloor a_t \rfloor + k}-x\right) \leq C_2e^{-\alpha x}.
\end{align*}
Equation $\eqref{switch}$ yields:
\begin{equation*}
\mathbb{P}\left(\overline{X}\left( t_{\lfloor a_t \rfloor + k}-x\right)\geq \lfloor a_t \rfloor + k\right) \leq C_2e^{-\alpha x}. 
\end{equation*}
By taking $x=t_{\lfloor a_t \rfloor + k}-t$, we obtain:
\begin{equation*}
\mathbb{P}\left(\overline{X}(t)-\lfloor a_t \rfloor \geq  k\right) \leq C_2e^{-\alpha \left(t_{\lfloor a_t \rfloor + k}-t\right)}. 
\end{equation*}
A straightforward computation shows that:
\begin{equation*}
t_{\lfloor a_t \rfloor + k}-t=\frac{1}{c^+}\left[k+\lfloor a_t \rfloor-a_t+\frac{3}{2\theta^+}\log(c^+)\right]+o_t(1),
\end{equation*}
and thus the first part of \eqref{distingxatk} follows. Since the proofs of the two parts of \eqref{distingxatk} are similar, we omit the second. Equation \eqref{asymptailH} shows that the distribution of the centred maximum has exponentially bounded asymptotics. Moreover, we precisely know the asymptotics of $\overline{\phi}$ by \cite{ChRoConnecting} and $\eqref{asympgauchedrmota}$: 
\begin{equation} 
1-\overline{\phi}(x) \underset{x \rightarrow + \infty}{\sim} xe^{-\theta x} \quad \mbox{ and } \quad \overline{\phi}(x)  \underset{x \rightarrow - \infty}{\sim} e^{\frac{ x}{c^+}}. \label{2asympphi}
\end{equation} 
Using these facts we will prove \eqref{firstpartcoroperio}. Let $\delta>0$, such that $\delta<\min\left(\theta,\frac{1}{c^+},\alpha^{\prime}\right)$. Let $f$ be a piecewise continuous function, such that $f(x)= \underset{x \rightarrow \pm \infty}{o}(e^{\delta |x|})$ and for $k \in \mathbb{Z}$, let $g_k $ and $p_k$ be defined by:
\begin{align} g_k(t)&= f(k-a_t)\; \mathbb{P}\left(\overline{X}(t) = k\right) \nonumber \\
&=  f(k-a_t)\left[ \overline{h}(k,t)-\overline{h}(k-1,t)\right], \quad \forall t>0 \label{defgk}
\end{align}
and
\begin{equation}p_k(s)= f(k-s)\left[ \overline{\phi}\left(k-s\right)-\overline{\phi}\left(k-1-s\right)\right], \quad \forall s \in \mathbb{R}.
\end{equation}
We will first show that $G_f$ and $P_f$ are well-defined. For $k \in \mathbb{Z}$ and $t \geq \eta$, we have:
 \begin{align}
 \mathbb{P}\left(\overline{X}(t) = k+\lfloor a_t \rfloor\right) &\leq  \mathbb{P}\left(\left|\overline{X}(t) -\lfloor a_t \rfloor \right|=|k|\right) \nonumber \\
&\leq  \mathbb{P}\left(\left|\overline{X}(t) -\lfloor a_t \rfloor \right| \geq |k|\right) \nonumber \\
&\leq  C_3e^{-\alpha^{\prime} |k|}, \label{ineqpxmaxegalkat}
 \end{align}
 where \eqref{ineqpxmaxegalkat} is a direct consequence of \eqref{asymptailH}. Furthermore, the assumptions on $f$ and \eqref{ineqpxmaxegalkat} imply that there exists $C_6>0$ such that:
\begin{equation}
g_{k+\lfloor a_t \rfloor}(t)  \leq C_6e^{-(\alpha^{\prime}-\delta)|k|}, \quad \forall k \in \mathbb{Z}, \forall t \geq \eta.\label{asmphandg}
\end{equation}
Similarly the asymptotics in $\eqref{2asympphi}$ and the assumptions on $f$ yield that there exists $(\delta^{\prime},C_7)\in {(\mathbb{R}^*)}^2$ such that:
\begin{equation}
 \quad p_{k+\lfloor a_t \rfloor}(a_t)  \leq C_7e^{-\delta^{\prime}|k|}, \quad \forall k \in \mathbb{Z}, \forall t \geq \eta. \label{asmphandp}
\end{equation}


We see that the sums $P_f$ and $G_f$ defined in \eqref{gf1} are invariant by translation by an integer. Thus:
\begin{equation}
G_f(t)=\sum_{k \in \mathbb{Z}} g_{k+\lfloor a_t \rfloor}(t) \quad \mbox{ and } \quad r_f(t)=P_f(a_t)=\sum_{k \in \mathbb{Z}} p_{k+\lfloor a_t \rfloor} (a_t). \label{GfandPFtranslted}
\end{equation}
The upper bounds $\eqref{asmphandg}$ and $\eqref{asmphandp}$ ensure that $G_f$ and $r_f$ are well-defined on $[\eta,+\infty)$, for each $\eta>0$. Therefore, $G_f$ and $r_f$ are well-defined on $\mathbb{R}^{*}_+$. We can also observe that if we further suppose that $f$ is continuous on $\mathbb{R}$, since $\eqref{asmphandg}$ and $\eqref{asmphandp}$ yields the normal convergence of the series  in \eqref{GfandPFtranslted}, $G_f$ and $r_f$ are continuous and by a simple change of variables it is also the case for $P_f$.  Let us now prove \eqref{firstpartcoroperio}. The upper bounds $\eqref{asmphandg}$ and $\eqref{asmphandp}$  yield that for $\epsilon>0$, there exist $k_0 \in \mathbb{Z}_{-}$ and $k_1 \in \mathbb{N}$ such that for all $t \geq \eta$:
\begin{equation}
\sum_{k=-\infty}^{k_0} g_{k+\lfloor a_t \rfloor}(t) + \sum_{k =k_1}^{+\infty} g_{k+\lfloor a_t \rfloor}(t)< \epsilon, \label{sumkminusinfh}
\end{equation}
and
\begin{equation}
\sum_{k =-\infty}^{k_0} p_{k+\lfloor a_t \rfloor}(a_t) + \sum_{k =k_1}^{+\infty} p_{k+\lfloor a_t \rfloor}(a_t)< \epsilon. \label{sumkminusinfp}
\end{equation}
Furthermore, since there is a finite number of integers between $k_0$ and $k_1$, Theorem $\ref{thconvergenceYule}$ ensures that :
\begin{equation}
\lim_{t \rightarrow +\infty} \sum_{k =k_0}^{k =k_1} |g_{k+\lfloor a_t \rfloor}(t)-p_{k+\lfloor a_t \rfloor}(a_t)|=0. \label{milieuth22}
\end{equation}
Equations \eqref{sumkminusinfh}, \eqref{sumkminusinfp} and \eqref{milieuth22} thus yield that for all $\epsilon>0$; there exists $t_0>0$ such that for all $t>t_0$
\begin{equation}
\sum_{k =-\infty}^{+\infty} |g_{k+\lfloor a_t \rfloor}(t)-p_{k+\lfloor a_t \rfloor}(a_t)|<3\epsilon,
\end{equation}
which is equivalent to \eqref{firstpartcoroperio}. \\

We have yet to deal with \eqref{secondpartcoroperio}. For $m \in \mathbb{N}$, let $f\in \mathcal{C}^{m}(\mathbb{R})$ with: 
\begin{equation}
f^{(i)}(x)=\underset{x \rightarrow \pm \infty}{o}(e^{\delta |x|}), \quad  \forall \; 0\leq i \leq m. \label{hypf}
\end{equation} 
The fact that $f \in \mathcal{C}^{m}(\mathbb{R})$ obviously implies that: $t \mapsto p_k(a_t)$ and $g_k$ belong to $\mathcal{C}^{m}(\mathbb{R_+^*})$. We want to prove that $G_f$ and $P_f$ belong to $\mathcal{C}^{m}(\mathbb{R_+^*})$. For that purpose, we will prove that,  for each $i \in \lbrace 1,...,m\rbrace$, the series $\left(\sum g^{(i)}_k\right)$ and $\left(\sum p^{(i)}_k\right)$ converge uniformly on every compact subset of $\mathbb{R}^*_+$. 

For $k \in \mathbb{Z}$, define $\tilde{f}_k$ by:
\begin{equation}
\tilde{f}_k(t)=f(k-a_t), \forall t>0.
\end{equation}
By an induction, we can show that for each $i \in \lbrace 1,...,m\rbrace$, there exists a sequence of $i-1$ polynomial functions  $(R_{j,i})_{j \in \lbrace 1,...,i-1 \rbrace}$ with $i$ variables such that:
\begin{equation}
\tilde{f}^{(i)}_k(t)=\sum_{j=1}^{i-1} R_{j,i}\left(a^{\prime}_t,...a^{(i)}_t\right)f^{(j)}(k-a_t)+(-a^{\prime}_t)^if^{(i)}(k-a_t) \label{sumftilde}
\end{equation}
and 
\begin{equation}
R_{j,i}\left(a^{\prime}_t,...a^{(i)}_t\right)=o_t(1), \forall j\in \lbrace 1,...,i-1 \rbrace.  \label{Rpetito}
\end{equation}
Equation  \eqref{Rpetito} is simply a consequence of the fact that $a^{\prime}_t=c^++o_t(1)$ and $a^{(j)}_t=o_t(1),$ $ \forall j \geq 2$.
We can also observe that the derivatives of $a$ are bounded on $[\eta,+\infty)$. Therefore the hypothesis $\eqref{hypf}$ on $f$ and Equation \eqref{sumftilde} give us the existence of $\tilde{C}_i>0$ such that we have:
\begin{equation}
\left|\tilde{f}^{(i)}_k(t)\right|\leq \tilde{C}_i e^{\delta |k-a_t|}, \quad \forall t\geq \eta. \label{majftilde} 
\end{equation}
We now want to give by induction an upper bound to $\partial^{(i)}_t \overline{h}(k,t)$. First suppose that $i=1$ and that $k\geq \lfloor a_t \rfloor$. Since $\partial_t \overline{h}(k,t)$ is negative Equation $\eqref{eqnpri}$ yields for $k\in \mathbb{Z}$ and $t\geq \eta$:
\begin{align}
\left|\partial_t \overline{h}(k,t)\right|&=\overline{h}(k,t)-\overline{h}^2(k-1,t) \nonumber \\
&\leq \overline{h}(k,t) \nonumber \\
&\leq \mathbb{P}\left(\overline{X}(t)\leq (k-\lfloor a_t \rfloor)+\lfloor a_t \rfloor\right) \nonumber\\
&\leq C_4e^{-\alpha^{\prime}(k-\lfloor a_t \rfloor)}, \label{fstmajorapartih}
\end{align}
where the last inequality comes from \eqref{distingxatk}. Similarly, if we now suppose that $k\leq \lfloor a_t \rfloor$, we get:
\begin{align}
\left|\partial_t \overline{h}(k,t)\right|&=1-\mathbb{P}\left(\overline{X}(t)\geq k\right)-\left(1-\mathbb{P}\left(\overline{X}(t)\geq k-1\right)\right)^2\nonumber \\
&\leq 2\mathbb{P}\left(\overline{X}(t)\geq k-1\right) \nonumber \\
&\leq 2e^{-\alpha^{\prime}}C_5e^{-\alpha^{\prime}(\lfloor a_t\rfloor-k )}, \label{sndmajorapartih}
\end{align}
where $C_5$ is defined in \eqref{distingxatk}. Therefore, Equations \eqref{fstmajorapartih} and \eqref{sndmajorapartih} and the fact that $a_t-\lfloor a_t \rfloor\leq 1$ yield that there exists $C_8>0$ such that:
\begin{equation}
\left|\partial_t \overline{h}(k,t)\right|\leq C_8 e^{-\alpha^{\prime}|k-a_t|}.\label{majpartialhbyC8}
\end{equation}
Fix $i\geq 2$ and suppose now that for each  $j \in \lbrace 1,...,i-1 \rbrace$, there exists $\tilde{C}^{\prime}_{j}>0$ such that for every $k\in \mathbb{Z}$ and $t\geq \eta$:
\begin{equation}\left|\partial^{(j)}_t \overline{h}(k,t) \right| \leq \tilde{C}^{\prime}_j e^{-\alpha^{\prime}|k- a_t |}. \label{majderiv}
\end{equation}
If we differentiate $i-1$ times with respect to $t$ Equation $\eqref{eqnpri}$ we get for $k \in \mathbb{Z}$ and $t>0$:
\begin{equation}
\partial^{(i)}_t \overline{h}(k,t)=\sum_{j=0}^{i-1} \binom{i-1}{j} \partial^{(j)}_t \overline{h}(k-1,t) \partial^{(i-j-1)}_t \overline{h}(k-1,t)-\partial^{(i-1)}_t \overline{h}(k,t). \label{deriveeiemeh}
\end{equation}  

Since $\left|\overline{h}\right|\leq 1$ and using the upper bounds \eqref{majderiv} into \eqref{deriveeiemeh}, we get that there exists $\tilde{C}^{\prime}_i>0$ such that \eqref{majderiv} holds with $j=i$. We thus have proved by induction that  \eqref{majderiv} holds for all $j\in \lbrace 1,...,m \rbrace$. 

Fix $i\in \lbrace 1,...,m \rbrace$. By differentiating $i$ times $g_k$ defined in \eqref{defgk}, we get:
\begin{equation}
g^{(i)}_k(t)=\sum_{j=0}^{i}\binom{i}{j} \tilde{f}^{(j)}_k(t)\left(\partial^{(i-j)}_t \overline{h}(k,t) -\partial^{(i-j)}_t \overline{h}(k-1,t)\right) , \quad \forall t>0. \label{gkderivativei}
\end{equation}
By introducing the upper bounds \eqref{majftilde} and \eqref{majderiv} into \eqref{gkderivativei}, we obtain that there exists $\tilde{C}^{(2)}_i>0$ such that:
\begin{equation}
\left|g^{(i)}_k(t)\right|\leq \tilde{C}^{(2)}_i e^{-(\alpha^{\prime}-\delta) |k-a_t|}, \quad \forall t \geq \eta. \label{majga}
\end{equation}
Since $a$ is smooth on $\mathbb{R}^*_+$, it is bounded on every compact subset of $\mathbb{R}^*_+$. Hence, \eqref{majga}   implies that the series $\left(\sum g^{(i)}_k\right)$ converges uniformly on every compact subset of $\mathbb{R}^*_+$ for each $i \in \lbrace 1,...,m \rbrace$. We have already proved that it is also the case for $i=0$. Therefore, $G_f \in \mathcal{C}^m(\mathbb{R}^*_+)$. With the same approach, we can show that $P_f \in \mathcal{C}^m(\mathbb{R})$. \\

Let us now prove \eqref{secondpartcoroperio}. Recall the definition of $r_f$ in \eqref{gf1}. As for $\tilde{f}$, a induction yields for each $i \in \lbrace 1,...,m\rbrace$, the existence of a sequence of $i-1$ polynomial functions $(\tilde{R}_{j,i})_{j \in \lbrace 1,...,i-1 \rbrace}$ with $i$ variables such that:
\begin{equation}
r^{(i)}_f(t)=a^{(i)}_tP^{(i)}_f(a_t)+\sum_{j=0}^{i-1}\tilde{R}_{j,i}\left(a^{\prime}_t,...,a^{(i-1)}_t\right)P^{(j)}_f(a_t) \label{tildeRandr}
\end{equation}
and
\begin{equation}
\tilde{R}_{j,i}\left(a^{\prime}_t,...,a^{(i-1)}_t\right)=o_t(1), \; \forall j \in \lbrace 1,...,i-1 \rbrace. \label{tildeRpetito}
\end{equation}
 Since the derivatives of $P_f$ are periodic, they are bounded. Therefore \eqref{tildeRandr} and \eqref{tildeRpetito} yield:
\begin{equation}
r^{(i)}_f(t)=a^{(i)}_tP^{(i)}_f(a_t)+o_t(1)={(c^+)}^iP^{(i)}_f(a_t)+o_t(1).\label{tildeRandr2}
\end{equation}
We recall that $P_f$ and $G_f$ are invariant by a translation by an integer. Therefore:
\begin{align}
r^{(i)}_f(t)-G^{(i)}_f(t)&=\sum_{k=-{\infty}}^{+\infty}\left[{(c^+)}^ip^{(i)}_k(a_t)-g^{(i)}_k(t) \right]+o_t(1)\nonumber \\
&=\sum_{k=-{\infty}}^{+\infty}\left[{(c^+)}^ip^{(i)}_{k+\lfloor a_t \rfloor}(a_t)-g^{(i)}_{k+\lfloor a_t \rfloor}(t)\right]+o_t(1). \label{thesumderiverG}
\end{align}
where
\begin{align}
{(c^+)}^ip^{(i)}_{k+\lfloor a_t \rfloor}(a_t)=\sum_{j=0}^{i}&\left[\binom{i}{j} f^{(j)}(k-\lbrace a_t\rbrace){(-c^+)}^{i} \times\right. \nonumber\\
&\left. \left(\overline{\phi}^{(i-j)}(k-\lbrace a_t\rbrace) - \overline{\phi}^{(i-j)}(k-1-\lbrace a_t\rbrace)\right)\right]
\end{align}
and
\begin{align}
g^{(i)}_{k+\lfloor a_t \rfloor}(t)=\sum_{j=0}^{i}&\left[\binom{i}{j}f^{(j)}(k-\lbrace a_t\rbrace)(-c^+)^{j}  \times \right. \nonumber\\
&\left.\left(\partial^{(i-j)}_t \overline{h}(k+\lfloor a_t \rfloor,t) -\partial^{(i-j)}_t \overline{h}(k+\lfloor a_t \rfloor-1,t)\right)\right]+o_t(1).
\end{align}
To prove the final part of the Theorem, we can again cut the sum in \eqref{thesumderiverG} in three terms and proceed as for the proof of \eqref{firstpartcoroperio}. We will simply show by induction that:
\begin{equation}
\lim_{t \rightarrow + \infty} \sup_{k \in \mathbb{Z}} \left|\partial_t^{(i)} \overline{h}(k+\lfloor a_t \rfloor, t)-(-c^+)^{i}\overline{\phi}^{(i)}(k -\{a_t\} )\right|=0. \label{convergencederiveiem}
\end{equation}
The base case $i=0$ is Theorem $\ref{thconvergenceYule}$. If we assume that the result holds for the rank $i-1$, we can rewrite $\eqref{deriveeiemeh}$:
\begin{align}
\frac{\partial^{(i)}_t \overline{h}(k,t)}{(-c^+)^{i-1}}=&\sum_{j=0}^{i-1} \binom{i-1}{j} \overline{\phi}^{(j)}(k-1-\{a_t\} )\overline{\phi}^{(j-i-1)}(k-1-\{a_t\} ) \nonumber \\
&-\overline{\phi}^{(i-1)}(k-\{a_t\} ) +o_t(1),
 \label{derivee2iemeh}
\end{align}  
where $o_t(1)$ is uniform in $k$. Differentiating $i-1$ times Equation $\eqref{eqx}$  and combining the result with $\eqref{derivee2iemeh}$, we obtain  $\eqref{convergencederiveiem}$. 

The end of the proof is identical as for the proof of \eqref{firstpartcoroperio}.
 \end{proof}
For the minimum the proof is the same except that one of the asymptotics of $\underline{\phi}$ $\eqref{2asympphi}$ is different from those of $\overline{\phi}$. Indeed, $\underline{\phi}$ is decreasing and thus $\eqref{eqx}$ implies that for all $x \in \mathbb{R}$, $\underline{\phi}(x)\leq \underline{\phi}^2(x-1)$. Consequently, there is $0<A<1$ such that: $$\underline{\phi}(x) \leq A^{2^x}, \quad \forall x >0.$$
\subsection{Proof of Corollary \ref{propuprim} }
We just prove Lemma \ref{lemmawt}. Indeed, if we suppose that Lemma \ref{lemmawt} holds, and if we define $f$ and $g$ by $f(x)=x$ and $g(x)=x^2$, Theorem \ref{coroperio} yields:
$$
\begin{cases}
\mathbb{E}\left(\overline{X}(t)-a_t\right)=G_f(t) \\
\mbox{Var}\left(\overline{X}(t)\right)= \mbox{Var}\left(\overline{X}(t)-a_t\right)=G_g(t)-G^2_f(t)\\
\mathbb{E}\left(\tilde{F}_t\right)=G^{\prime}_f(t)+a^{\prime}_t=G^{\prime}_f(t)+c^++o_t(1) \\
\end{cases}
$$
and Corollary \ref{propuprim} immediately follows.
\begin{proof}[Proof of Lemma \ref{lemmawt}]
Let $t\geq 0$ and $h>0$. We first recall that Theorem \ref{coroperio} applied to the identity ensures that $w:t \mapsto \mathbb{E}\left(\overline{X}(t)\right)$ is well-defined and differentiable. We will show that $w^{\prime}(t)=\mathbb{E}\left(\tilde{F}_t\right)$. For this purpose, we divide $\overline{X}(t+h)-\overline{X}(t)$ into three parts depending on the number $J_{t,h}$ of jump of $(N_t)$ between $t$ and $t+h$. When a particle dies, it gives birth to 2 particles and thus $J_{t,h}=N_{t+h}-N_t$.  

First consider $\nu_0(t,h):=\left(\overline{X}(t+h)-\overline{X}(t)\right)\mathbf{1}_{\lbrace N_{t+h}-N_t=0 \rbrace}.$ If there is no division between $t$ and $t+h$, $\overline{X}(t+h)=\overline{X}(t)$. Therefore, $\nu_0(t,h)=0$. 

Now set $\nu_1(t,h):=\left(\overline{X}(t+h)-\overline{X}(t)\right)\mathbf{1}_{\lbrace N_{t+h}-N_t=1 \rbrace}.$ By construction of the Yule tree process, when there is a division, the particle which splits is chosen uniformly and independently of the number of particles. Furthermore, $\overline{X}$ increases by 1 after a split if and only the split occurs for a particle situated at the maximal position. Therefore,
\begin{equation}
\mathbb{E}(\nu_1(t,h)|\mathcal{F}_t)=\frac{\tilde{F}_t}{N_t}\mathbb{P}\left(N_{t+h}-N_t=1|\mathcal{F}_t\right).
\end{equation}
Since $(N_t)$ is a pure birth process, we have:
\begin{equation}
\mathbb{P}\left(N_{t+h}-N_t=1|\mathcal{F}_t\right)=N_t h +o_h(h), \mbox{  uniformly for all } t. 
\end{equation}
Consequently,
\begin{equation}
\mathbb{E}(\nu_1(t,h))=h\mathbb{E}(\tilde{F}_t)+o_h(h).
\end{equation}
Finally, consider $\nu_2(t,h)=\left(\overline{X}(t+h)-\overline{X}(t)\right)\mathbf{1}_{\lbrace N_{t+h}-N_t\geq 2 \rbrace}$. Since at each split $(N_t)$ increases by one and $\overline{X}$ increases at most by one, we have
\begin{equation}
\nu_2(t,h)\leq\left(N_{t+h}-N_t\right)\mathbf{1}_{\lbrace N_{t+h}-N_t\geq 2 \rbrace} \\
\end{equation}
By replacing $\mathbf{1}_{\lbrace N_{t+h}-N_t\geq 2 \rbrace} $ by $1-\mathbf{1}_{\lbrace N_{t+h}-N_t=1 \rbrace} -\mathbf{1}_{\lbrace N_{t+h}-N_t=0 \rbrace} $, we obtain
\begin{equation}
\left(N_{t+h}-N_t\right)\mathbf{1}_{\lbrace N_{t+h}-N_t\geq 2 \rbrace}= N_{t+h}-N_t-\mathbf{1}_{\lbrace N_{t+h}-N_t=1 \rbrace}.\label{G2NTh}\\
\end{equation}
By taking the expectation of the terms of Equation \eqref{G2NTh}, we get:
\begin{align*}
\mathbb{E}\left(\left(N_{t+h}-N_t\right)\mathbf{1}_{\lbrace N_{t+h}-N_t\geq 2 \rbrace}\right)&= e^{t+h}-e^t-he^t+o_h(h)\\
&= o_h(h).
\end{align*}
Therefore, 
\begin{equation}
\mathbb{E}(\nu_2(t,h))=o_h(h).
\end{equation}
By grouping $\nu_0$, $\nu_1$ and $\nu_2$, we get:
\begin{equation}
\mathbb{E}(\overline{X}(t+h)-\overline{X}(t))=h\mathbb{E}(\tilde{F}_t)+o_h(h),
\end{equation}
which concludes the proof.
\end{proof}
\section{Application to the binary search tree} \label{sectbinasear}
 We recall that by  \eqref{embed1} the binary search tree is embedded into the Yule tree. Consequently, if $H_n$ is the height of a random binary search tree and if $\tau_n$ is defined as in \eqref{tau}, we have:
 \begin{equation}
 H_n = \overline{X}(\tau_n). \label{equimax}
 \end{equation}
Drawing our inspiration from Lalley and Selke \cite{lalley1987conditional} and with the help of Theorem \ref{thconvergenceYule}, we will prove  Theorem $\ref{thconvergencearbre}$.
\begin{proof}[Proof of Theorem $\ref{thconvergencearbre}$]
For $x\in \mathbb{R}$, recall that $\{x\} :=x-\lfloor x \rfloor$. Our proof is divided into two steps. \\ \\
Step 1: First, we will show that: \\ \\
$\exists K>0: \forall \epsilon>0 \; \exists  s_0>0 : \forall s>s_0 \; \exists t_s>0 : \forall t>t_s: $
\begin{align}
\bigg \vert &\mathbb{P}\left(\overline{X}\left(t+s-\log(W_{s})\right) \leq \lfloor a_{t+s}+x \rfloor \right)  \nonumber \\
&\qquad \qquad - \mathbb{E}\left(\exp\left(-Ke^{-\theta (\lfloor a_{t+s}+x \rfloor -  a_{t+s})} \partial Z^+_{\infty}\right)\right)\bigg\vert < \frac{\epsilon}{3}, \label{Binstep1}
\end{align}
where $W_s=W_s(0)$. By Markov property, we have for $t,s>0$ and $x \in \mathbb{R}$ that:
\begin{align}
&\mathbb{P}\left(\overline{X}\left(t+s-\log(W_s)\right) \leq a_{t+s}+x  \big\vert \mathcal{F}_s \right) \nonumber \\
&=\prod_{u \in \mathcal{N}_s} \mathbb{P}\left(\overline{X}^u\left(t-\log(W_s)\right) \leq a_{t+s}+x -X_u(s) \right) \nonumber \\
&=\prod_{u \in \mathcal{N}_s} \mathbb{P}\left(\overline{X}^u\left(t-\log(W_s)\right) \leq a_{t-\log{W_s}} + a_{t+s}+x -X_u(s) -a_{t-\log W_s}\right),
\end{align}
where the processes $\overline{X}^u$, defined by $\overline{X}^u(t)=\max \lbrace X_v(t+s), v \in \mathcal{N}_{t+s},u<v \rbrace$, $\forall t>0$, $\forall u \in \mathcal{N}_s$ are independent of $\mathcal{F}_s$ and identically distributed.
For a fixed $s$, we have by Theorem $\ref{thconvergenceYule}$ that almost surely:
\begin{align}
&\lim_{t \rightarrow +\infty} \bigg \vert \mathbb{P}\left(\overline{X}\left(t+s-\log(W_s)\right) \leq a_{t+s}+x  \big\vert \mathcal{F}_s \right) \nonumber  \\
&\qquad \quad -\prod_{u \in \mathcal{N}_s} \phi_K( \lfloor a_{t+s} +x \rfloor -a_{t-\log W_s}-X_u(s))\bigg \vert =0,
\end{align}
where $\phi_K(x)=\mathbb{E}\left(\exp\left(-Ke^{-\theta^+ x} \partial W_{\infty} (\theta^+)\right)\right)=\overline{\phi}(x)$. Observe that the argument of $\phi_K$ can be rewritten as:
\begin{align*}
\lfloor a_{t+s} +x \rfloor -a_{t-\log W_s}-X_u(s) &= a_{t+s}+x -\{a_{t+s}+x\}-a_{t-\log(W_s)}-X_u(s) \\
&=c^+s+c^+\log W_s +R_{x,s,t}+S_{s,t}-X_u(s),
\end{align*}
where 
$$R_{x,s,t}=x-\{a_{t+s}+x\}\quad \mbox{ and } \quad S_{s,t}=\frac{3 \log(\frac{t-\log W_s}{t+s})}{2 \theta^+}.$$ Moreover, $S_{s,t}$ goes to $0$ when $t$ goes to infinity. $\phi_K$ is continuous and bounded by $1$, and therefore by dominated convergence:
\begin{align}
&\lim_{t \rightarrow +\infty} \bigg \vert \mathbb{P}\left(\overline{X}\left(t+s-\log(W_s)\right) \leq a_{t+s}+x  \right) \nonumber  \\
&\qquad \quad -\left. \mathbb{E}\left(\prod_{u \in \mathcal{N}_s} \phi_K\left( c^+s+c^+\log W_s +R_{x,s,t}-X_u(s)\right) \right) \right| =0. \label{BinStep11}
\end{align}
We will now rewrite the product of the terms from $\eqref{BinStep11}$:
\begin{align}
&\prod_{u \in \mathcal{N}_s} \phi_K(c^+s+c^+\log W_s +R_{x,s,t}-X_u(s) ) \\
=&\prod_{u \in \mathcal{N}_s} \mathbb{E}\left. \left(\exp\left(-Ke^{-\theta\left(c^+s+c^+\log W_s +R_{x,s,t}-X_u(s)\right)}\partial W_{\infty,u}(\theta^+)\right)\right|\mathcal{F}_s\right) \\
=& \quad \mathbb{E}\left.\left(\exp\left(-Ke^{-\theta\left(c^+\log W_s +R_{x,s,t}\right)}\sum_{u \in \mathcal{N}_s} e^{\theta\left(X_u(s)-c^+s\right)}\partial W_{\infty,u}(\theta^+)\right)\right| \mathcal{F}_s\right),
\end{align}
where $\partial W_{\infty,u}(\theta^+)$ are independent copies of $\partial W_{\infty}(\theta^+)$ and independent of $\mathcal{F}_s$. Moreover we know (\cite{ChRoConnecting} and \eqref{embedmartingalebstandyule1}) that: 
\begin{equation}
\sum_{u \in \mathcal{N}_s} e^{\theta\left(X_u(s)-c^+s\right)}\partial W_{\infty,u}(\theta^+)=\partial W_{\infty}(\theta^+)=\frac{W_{\infty}(0)^{c^++1}}{\Gamma(c^+)}\partial Z_{\infty}^+.
\end{equation}
The definition of $c^+$ means that $c^+\theta^+=c^++1$ and thus:
\begin{align}
&\prod_{u \in \mathcal{N}_s} \phi_K(c^+s-X_u(s)+c^+\log W_s +R_{x,s,t} ) \nonumber \\
=& \quad \left. \mathbb{E}\left(\exp\left(-Ke^{-\theta^+\left(c^+\log W_s +R_{x,s,t}\right)}\frac{W_{\infty}(0)^{c^++1}}{\Gamma(c^+)}\partial Z_{\infty}^+\right)\right| \mathcal{F}_s\right)  \nonumber \\
=& \quad \left. \mathbb{E}\left(\exp\left(-K^{\prime}e^{-\theta^+ R_{x,s,t}}\left(\frac{W_{\infty}(0)}{W_s(0)}\right)^{c^++1}\partial Z_{\infty}^+\right)\right| \mathcal{F}_s\right), \label{Kprime}
\end{align}
with $K^{\prime}=\frac{K}{\Gamma(c^+)}$.
Let $\mu_{K^{\prime}}(R,s)=\mathbb{E}\left(\exp\left(-K^{\prime}e^{-\theta^+ R}\left(\frac{W_{\infty}(0)}{W_s(0)}\right)^{c^++1}\partial Z_{\infty}^+\right)\right)$.
By dominated convergence:
$$\lim_{s \rightarrow + \infty} \mu_{K^{\prime}}(R,s) = \psi(R).$$
Moreover for fixed $s$, $\mu(.,s)$ is increasing and
$$\lim_{x \rightarrow + \infty} \mu(x,s)=1 \quad \mbox{ and } \quad \lim_{x \rightarrow - \infty} \mu(x,s)=0.$$ 
Thus by Theorem \ref{polyatheorem}:
$$\lim_{s \rightarrow + \infty} \sup_{R \in \mathbb{R}} |\mu_{K^{\prime}}(R,s) - \psi(R)|=0.$$
We then have:
\begin{equation}
\lim_{s \rightarrow + \infty} \sup_{t > 0} \bigg \vert \mathbb{E}\left(\prod_{u \in \mathcal{N}_s} \phi_K(c^+s-X_u(s)+c^+\log W_s +R_{x,s,t} )\right)-\psi_{K^{\prime}}(R_{x,s,t}) \bigg \vert =0.\label{BinStep12}
\end{equation}
So from Equations $\eqref{BinStep11}$ and $\eqref{BinStep12}$ we deduce $\eqref{Binstep1}$. \\

Step 2: Noticing that $\log(n)-\log(W_{\tau_n}(0))=\tau_n$, we have by $\eqref{equimax}$:
\begin{align}
\mathbb{P}(H_n \leq \lfloor a_{\log(n)}+x \rfloor )&=\mathbb{P}(\overline{X}(\tau_n)\leq  \lfloor a_{\log(n)}+x \rfloor ) \nonumber\\
&=\mathbb{P}(\overline{X}(\log(n)-\log(W_{\tau_n}(0)))\leq  \lfloor a_{\log(n)}+x \rfloor ). \label{Binstep21}
\end{align}
For every $K>0$, $\psi_K$ is uniformly continuous and hence for $\epsilon>0$ there is $\eta>0$ such that for all $(x,y) \in \mathbb{R}$ such that $|x-y|<\eta$, 
\begin{equation}
|\psi_K(x)-\psi_K(y)|<\frac{\epsilon}{3}. \label{continupsi}
\end{equation}
For $\eta^{\prime}>0$, we may choose $s_1>s_0$ (where $s_0$ is defined in $\eqref{Binstep1}$)  such that for every $s>s_1$ : 
\begin{equation}
\mathbb{P}\left(|\log(W_s(0))-\log(W_{\infty}(0))|>\frac{\eta^{\prime}}{2}\right)<\frac{\epsilon}{6}. \label{approxWs}
\end{equation}
 by the almost sure convergence of $W_s(0)$ to a positive random variable. In the same way, there exists $n_1 \in \mathbb{N}$ such that: $\log(n_1)>t_{s_1}+s_1$ and such that for all $n \geq n_1$: 
 \begin{equation}
 \mathbb{P}\left(|\log(W_{\tau_n}(0))-\log(W_{\infty}(0))|>\frac{\eta^{\prime}}{2}\right)<\frac{\epsilon}{6}. \label{WtnWinf}
 \end{equation}
 Introducing the approximations $\eqref{approxWs}$ and $\eqref{WtnWinf}$ in Equation $\eqref{Binstep21}$ and using the monotonicity of $\psi_{K^{\prime}}$ (where $K^{\prime}$ is defined in $\eqref{Kprime}$) we obtain for $n \geq n_1$:
 \begin{align}
 \mathbb{P}(H_n \leq \lfloor a_{\log(n)}+x \rfloor ) &\leq \mathbb{P}(\overline{X}(\log(n)-\log(W_{s_1}(0))-\eta^{\prime})\leq  \lfloor a_{\log(n)}+x \rfloor ) +\frac{\epsilon}{3} \nonumber\\
&\leq \psi_{K^{\prime}}( \lfloor a_{\log(n)}+x \rfloor -a_{\log(n)-\eta^{\prime}})+\frac{2\epsilon}{3}, \label{Step1interv}
 \end{align}
 where Equation $\eqref{Step1interv}$ is a consequence of $\eqref{Binstep1}$. So, taking $\eta^{\prime}>0$ such that for all $n\geq n_1$, $|a_{\log(n)-\eta^{\prime}}-a_{\log(n)}|<\eta$, we obtain:
 \begin{align}
 \mathbb{P}(H_n \leq \lfloor a_{\log(n)}+x \rfloor ) &\leq \psi_{K^{\prime}}( \lfloor a_{\log(n)}+x \rfloor -a_{\log(n)}+\eta)+\frac{2\epsilon}{3} \nonumber \\
&\leq \psi_{K^{\prime}}( \lfloor a_{\log(n)}+x \rfloor -a_{\log(n)})+\epsilon,  
 \end{align}
 by $\eqref{continupsi}$. Similarly, we obtain:
 \begin{equation}
\psi_{K^{\prime}}( \lfloor a_{\log(n)}+x \rfloor -a_{\log(n)})-\epsilon \leq \mathbb{P}(H_n \leq \lfloor a_{\log(n)}+x \rfloor ), 
 \end{equation}
 and then for $x \in \mathbb{R}$,
 $$\lim_{n \rightarrow +\infty} | \mathbb{P}(H_n \leq \lfloor a_{\log(n)}+x \rfloor ) -\psi_{K^{\prime}}( \lfloor a_{\log(n)}+x \rfloor -a_{\log(n)})|=0.$$
 By the same arguments as in the proof of Theorem $\ref{thconvergenceYule}$, we have, in fact, that:
 \begin{equation}
\lim_{n \rightarrow +\infty} \sup_{x \in \mathbb{R}}| \mathbb{P}(H_n \leq \lfloor a_{\log(n)}+x \rfloor ) -\psi_{K^{\prime}}( \lfloor a_{\log(n)}+x \rfloor -a_{\log(n)})|=0.
\end{equation}
\end{proof}
\begin{appendices} 
\section{Travelling-waves and martingales}
\subsection{Travelling-waves and martingales of the Yule branching random walk}
It is well-known that travelling waves and some martingales play a key role in the study of the extremal particles of the branching random walk. In this appendix we define the relevant objects in our context and recall the pertinent results.

 \begin{thma} \label{thexistenceunicite}
Equation $\eqref{eqx}$ has monotone and bounded travelling-wave solutions at speed $c$, in $\mathcal{C}^1(\mathbb{R})$ if and only if $c \leq c^-$ or  $c \geq c^+$. Moreover, uniqueness holds for each such $c \neq 0$ (up to an additive constant in the argument). For $c=0$, the whole set of solutions is:
\begin{equation}
\lbrace x \mapsto P(x)^{2^x} | \; P \mbox{ is 1-periodic and nonnegative}\rbrace. \label{twc0}
\end{equation}
\end{thma}
It is easy to see that we can find more than one monotone bounded solution in the set $\eqref{twc0}$. \\

\begin{Rem}
Let us point out that this result is stated in a slightly different form in \cite{ChRoConnecting}. Indeed, Chauvin and Rouault proved uniqueness for this equation up to a decreasing change of variables and in the class of Laplace transforms. The key to their proof was an application of Liu $\cite{liu1998fixed,liu2000generalized}$. More recently, Alsmeyer, Biggins and Meiners $\cite{alsmeyer2012functional}$ showed that there is  uniqueness among decreasing functions in $[0,1]$. 
We observe that the set of monotone bounded solutions of $\eqref{eqx}$ is in fact the set of monotone solutions in $[0,1]$. The proof that no travelling wave exist for $c \in (c^-,c^+)$ can be lifted as is from Harris\cite{harris1999travelling} for the branching Brownian motion. 
\end{Rem}
Since it is useful for our proofs, let us mention that the method  used in \cite{ChRoConnecting} also allows us to determine the left tail of the travelling-waves $\phi$ with speed $c\geq c^+$ of $\eqref{eqx}$:
\begin{equation}
\phi(x) \underset{x \rightarrow -\infty}{\sim}  e^{\frac{ x}{c}}. \label{asympgauchedrmota}
\end{equation} 
We now define the derivative martingale and the additive martingale. For that purpose, let us consider the natural filtration $(\mathcal{F}_t)_{t \geq 0}$ defined by $\mathcal{F}_t:= \sigma \lbrace X_u(s), u\in \mathcal{T}^c_s, \; 0 \leq s \leq t \rbrace, \; \forall t \geq 0$. We also fix sgn$(x)=1$, when $x\geq 0$ and sgn$(x)=-1$ when $x<0$.
 \begin{thmd}
 For $\theta \in \mathbb{R}$, the process $(W_t(\theta))_{t\geq 0}$ defined by
\begin{equation}
 W_t(\theta)=\sum_{u \in \mathcal{N}_t} e^{\theta(X_u(t)-c_{\theta} t)} \label{martaddi}
 \end{equation}
 is a $(\mathcal{F}_t)$-martingale called the additive martingale. Moreover, this martingale converges to an almost surely positive random variable $ W_{\infty}(\theta)$ when $\theta \in (\theta^-,\theta^+)$. In particular, the limit of the martingale is an exponential random variable with parameter 1 when $\theta=0$. \\
Similarly, the process $(\partial W_t(\theta))_{t\geq 0}$ defined by
 \begin{equation}
\partial W_t(\theta)=\mathrm{sgn}(\theta)\sum_{u \in \mathcal{N}_t}(2e^{\theta} t-X_u(t))e^{\theta(X_u(t)-c_{\theta} t)} \label{derivmart}
\end{equation}
 is a $(\mathcal{F}_t)$-martingale called the derivative martingale. Moreover, for $\theta \in \lbrace \theta^-,\theta^+\rbrace$, the limit $\partial W_{\infty}(\theta)$ of the derivative martingale exists and is almost surely positive.
\end{thmd}
Note that the derivative martingale of the Yule tree \eqref{derivmart} is the derivative with respect to $\theta$ (up to a change of sign for $\theta<0$) of the additive martingale. We now recall the link between travelling waves and these martingales.
\begin{thme}
 For $\theta \in (\theta^-,\theta^+)$, the travelling-wave at speed $c_{\theta}$ has the following representation:
 \begin{equation}
\phi_{K,\theta}(x)=\mathbb{E}\left(\exp\left({-Ke^{-\theta x} W_{\infty}(\theta)}\right)\right),\label{travnorm}
\end{equation}
where $K>0$ fixes the choice of the travelling wave. \\
For $\theta \in \lbrace\theta^-,\theta^+ \rbrace$, the travelling-wave at speed $c_{\theta}$ has the following representation:
\begin{equation}
\phi_{K,\theta}(x)=\mathbb{E}\left(\exp\left({-Ke^{-\theta x} \partial W_{\infty}(\theta)}\right)\right), \label{travcri}
\end{equation}
where $K>0$ fixes the choice of the travelling wave.
\end{thme}
\subsection{Martingales of the Yule generation process}
We now define the additive and derivative martingales for the Yule generation process. First consider the filtration $(\mathcal{F}^{GEN}_n)$ defined by $\mathcal{F}^{GEN}_n=\sigma \lbrace T_u(k),u \in \mathcal{M}_k, k\leq n \rbrace$. For $\theta \in \mathbb{R}$, let $(W^{GEN}_n(\theta))$ be defined by:
\begin{equation}
 W^{GEN}_n(\theta)=e^{\theta n}\sum_{u \in \mathcal{M}_n} e^{-\theta c_{\theta} T_u(n)}, \; \forall n \in \mathbb{N}. \label{martaddiyulegen}
 \end{equation}
 Similarly, define $(\partial W^{GEN}_n(\theta))$ by:
 \begin{equation}
 \partial W^{GEN}_n(\theta)=\mathrm{sgn}(\theta)e^{\theta n}\sum_{u \in \mathcal{M}_n}(2e^{\theta} T_u(n)-n) e^{-\theta c_{\theta} T_u(n)}, \; \forall n \in \mathbb{N}. \label{martderiyulegen}
 \end{equation}
 The following theorem groups together several results from \cite{ChRoConnecting}. We recall that $W_{\infty}(\theta)$ and $\partial W_{\infty}(\theta)$ are the limits of the additive and derivative martingales defined in \eqref{martaddi} and in \eqref{derivmart}.

 \begin{thma}
For $\theta \in \mathbb{R}$, the process $(W^{GEN}_n(\theta))$ is a $(\mathcal{F}^{GEN}_n)$-martingale called the additive martingale of the Yule generation process.  Similarly,  $(\partial W^{GEN}_n(\theta))$ is a $(\mathcal{F}^{GEN}_n)$-martingale called the derivative martingale of the Yule generation process.\\
Furthermore, for $\theta \in (\theta^-,\theta^+)$:
\begin{equation}
\lim_{n \rightarrow + \infty} W^{GEN}_n(\theta) = W_{\infty}(\theta) \; \mbox{a.s.} \label{convadditivgentotime}
\end{equation}
 and for $\theta \in \lbrace \theta^-,\theta^+\rbrace$:
 \begin{equation}
\lim_{n \rightarrow + \infty} \partial W^{GEN}_n(\theta) =\partial W_{\infty}(\theta) \; \mbox{a.s.} \quad .\label{convderivgentotime}
\end{equation}
 \end{thma}
 For the sake of brevity, we refrain from mentioning the idea of stopping lines in this article, we simply mention that the additive and the derivative martingales of the Yule generation process are the additive and the derivative martingales of the Yule branching random walk stopped on the sequence of stopping lines $\mathcal{M}_n$ which partially explains Equations \eqref{convadditivgentotime} and \eqref{convderivgentotime}. For more references about stopping lines see for instance \cite{chauvin1991product}.
 \subsection{Martingales of the binary search tree}
 Let us first define the additive martingale and the derivative martingale for the binary search tree. Let $\mathcal{L}_n$ be the set of leaves of $\mathcal{T}_n$ and define $\mathcal{F}_n^{BST}=\sigma(\lbrace u \in \mathcal{T}_i \rbrace , i \leq n, u \in \mathcal{T})$ where $\mathcal{T}$ is the complete binary tree. For $u \in \mathcal{L}_n$, we denote by $|u|$ the height of $u$.
For $z \in \mathbb{R}\setminus -\frac{\mathbb{N}}{2}$, we fix:
\begin{equation}
  U_0(z)=1 \quad \mbox{ and } \quad U_n(z)=\prod^{n-1}_{k=0} \frac{k+2z}{k+1}.
  \end{equation}
 Then, the process $(Z_n(z))$ introduced by Jabbour $\cite{jabbour2001martingales}$ and defined by
\begin{equation}
 Z_n(z)=\frac{1}{U_n(z)}\sum_{u \in \mathcal{L}_n} z^{|u|},
\end{equation}
is a $(\mathcal{F}_n^{BST})$-martingale, which we will call the additive martingale of the binary search tree.  The derivative martingale of the binary search tree is then simply defined by
\begin{equation}
\partial Z_n(z)=\frac{d Z_n(z)}{dz}. \label{martderiveearbrebinaire}
\end{equation}
The following theorem due to Chauvin and Rouault illuminates the connection between the Yule process and the derivative martingale $\eqref{martderiveearbrebinaire}$.
\begin{thma}
The martingale ($\partial Z_n(z)$) converges as $n \rightarrow + \infty$ to a positive random variable $\partial Z_{\infty}^+$ for $z=e^{\theta^+}$ (resp. to $\partial Z_{\infty}^-$ for $z=e^{\theta^-}$).
Moreover, by embedding the binary search tree into the Yule process, we have: 
\begin{equation}
 \partial W_{\infty}(\theta^+) =\frac{e^{-\theta^+}W_{\infty}(0)^{c^++1}}{\Gamma(c^+)}\partial Z_{\infty}^+ \label{embedmartingalebstandyule1}
 \end{equation}
 and:
 \begin{equation}
 \partial W_{\infty}(\theta^-) =\frac{e^{-\theta^-}W_{\infty}(0)^{c^-+1}}{\Gamma(c^-)}\partial Z_{\infty}^-,
 \end{equation}
 where the limit of the additive martingale $W_{\infty}(0)$ has a random exponential 1 law and is independent of $\partial Z_{\infty}^+$ and $\partial Z_{\infty}^+$.
\end{thma}
\section{Results on the maximum of a branching random walk}\label{sectionresultsonmaxim}
Since Addario-Berry and Reed's result \cite{addario2009minima} and Aïdékon's result \cite{aidekon2013convergence} are central to our argument, we here state them in more  details than in Section \ref{sectprevresult}. 
\subsection{Addario-Berry and Reed's result}\label{Subsecaddario}
 Consider a branching random walk defined as follows. A particle is at $0$ at time $0$. This particle dies at time $1$ and give birth to a random number of particles $B \in \mathbb{N}$ whose displacements are independent copies of a random variable $Y$. Then, for each $n\in \mathbb{N}$, each particle $u$ of the $n$th generation gives birth to $B_u$ particles, where $B_u$ is an independent copy of $B$, and the displacement of each new particle with respect to its parent is a copy of $Y$ independent of the others.  We also define $\Lambda$ by $\Lambda(\lambda):=\log\left(\mathbb{E}\left(e^{\lambda Y}\right)\right)$ and $ \overset{\circ}{D_{\Lambda}}$ the interior of the set of value for which $\Lambda(\lambda)$ is finite. Finally, we fix $M_n$ the minimum of the branching random walk and $\mathcal{S}$ the survival event. 
\begin{thmg}
Consider a branching random walk satisfying the following assumptions:
\begin{enumerate}
\item $\mathbb{E}(B)>1$,
\item there exists $d \geq 2$, such that $\mathbb{P}\left(B\leq d\right)=1$,
\item there exists $\lambda_0 > 0$, such that $\mathbb{E}\left(e^{\lambda_0 Y}\right)<+\infty$,
\item there exists $\lambda_1 \in \overset{\circ}{D_{\Lambda}}\cap \mathbb{R}^{*}_{-}$, such that $\lambda_1\Lambda^{\prime}(\lambda_1)-\Lambda(\lambda_1)=\log\left(\mathbb{E}(B)\right)$.
\end{enumerate}
Then 
\begin{equation}
\mathbb{E}(M_n|\mathcal{S})=\Lambda(\lambda_1) n-\frac{3}{2\lambda_1}\log n+O_n(1).
\end{equation}
Furthermore, there exist $C,\delta>0$ such that:
\begin{equation}\label{Addariobis}
\mathbb{P}(|M_n-\mathbb{E}(M_n|\mathcal{S})|\geq x|\mathcal{S})\leq Ce^{-\delta x}, \; \forall x \geq 0. 
\end{equation}
\end{thmg}
\subsection{Aïdékon's result}\label{SubsecAid}
The class of branching random walks in Aïdékon's article is slightly different. As before, a particle is at $0$ at time $0$. This particle dies at time $1$ and give birth this time to a non-lattice point process $\mathcal{L}$. Then, for each $n\in \mathbb{N}$, the particles of generation $n$ give birth to independent copies of $\mathcal{L}$, translated to their position. We call $\mathbb{T}$ the genealogical tree of the process and for each $u \in \mathbb{T}$, we denote by $|u|$ its generation and by $V(u)$ its position on the real line. We conserve the other notations of Addario-Berry and Reed's result. 

Aïdékon considers the boundary case, which is quite general after some renormalizations, and which corresponds to the following assumptions:
\begin{equation}\label{assumpaidekon1}
 \mathbb{E}\left(\sum_{|u|=1} 1 \right)>1, \;  \mathbb{E}\left(\sum_{|u|=1} e^{-V(u)} \right)=1,  \; \mathbb{E}\left(\sum_{|u|=1} V(u)e^{-V(u)} \right)=0.
 \end{equation}
 Furthermore, if we set $E_1:=\sum_{|u|=1} e^{-V(u)}$ and $E_2:=\sum_{|u|=1} V(u)_+e^{-V(u)}$, where $y_+=\max(0,y)$, the suppositions that:
 \begin{equation}\label{assumpaidekon2}
 \mathbb{E}\left(\sum_{|u|=1} V(u)^2e^{-V(u)}\right)<+\infty
 \end{equation}
and 
\begin{equation}\label{assumpaidekon3}
 \mathbb{E}\left(E_1\left(\log_+ E_1 \right)^2\right)<\infty,  \;  \mathbb{E}\left( E_2 \log_+ E_2 \right)<\infty
 \end{equation} 
 are made. We now introduce the derivative martingale:
\begin{equation}
 D_n:=\sum_{|u|=n}V(u)e^{-V(u)}, \; \forall n \in \mathbb{N}. \label{derivativeaidekon}
 \end{equation} 
It is know \cite{biggins2004measure} that $(D_n)$ converges to a positive limit $D_{\infty}$.
Finally, still writing  $M_n$ for the minimum of the branching random walk. Aïdékon shows under Assumptions \eqref{assumpaidekon1}, \eqref{assumpaidekon2} and \eqref{assumpaidekon3} the following theorem:
\begin{thmb}
There exists $C>0$ such that:
\begin{equation}
\lim_{n \rightarrow +\infty} \mathbb{P}\left(M_n \geq \frac{3}{2}\log(n)+x\right)=\mathbb{E}\left(e^{-Ce^xD_{\infty}}\right). \label{eqthAidekon}
\end{equation}
\end{thmb}
\section{A theorem of uniform convergence}
\begin{theorem}\label{polyatheorem}
Let $(f_n)$ be a sequence of non-decreasing functions from $\mathbb{R}$ to $[0,1]$, such that:
$$\lim_{x \rightarrow - \infty} f_n(x)=0 \mbox{ and } \lim_{x \rightarrow + \infty} f_n(x)=1.$$
If $(f_n)$ converges pointwise to a continuous function $f$ such that for all $n\in \mathbb{N}$:
$$\lim_{x \rightarrow - \infty} f(x)=0 \mbox{ and } \lim_{x \rightarrow + \infty} f(x)=1,$$ then the convergence is uniform.
\end{theorem}
This theorem is a simple extension of the following result (Problem 127 of \cite{polya1972g}).
\begin{theorem}\label{polyatheoremini}
Fix $a,b \in \mathbb{R}$ such that $a<b$. Let $(f_n)$ be a sequence of non-decreasing functions from $[a,b]$ to $\mathbb{R}$.
If $(f_n)$ converges pointwise to a continuous function $f$, then the convergence is uniform.
\end{theorem}
\begin{proof}[Proof of Theorem \ref{polyatheorem}]
Let $(f_n)$ and $f$ be defined as in the statement of Theorem \ref{polyatheorem}. Fix $\epsilon>0$. By hypothesis, there exists $x_0 \in \mathbb{R}$ such that $0\leq f(x_0)<\epsilon/6$. By pointwise convergence of $(f_n)$ to $f$ there exists $n_0$ such that for all $n\geq n_0$:
\begin{equation} \label{proofthpoyla1}
|f_n(x_0)-f(x_0)|<\epsilon/6.
\end{equation}
The choice of $x_0$ and Equation \eqref{proofthpoyla1} imply that for all $n \geq n_0$,
\begin{equation} 
0\leq f_n(x_0)<\epsilon/3.
\end{equation}
Since $(f_n)$ is a sequence of non-decreasing functions and since the limit of such a sequence is itself non-decreasing, we get that for all $n\geq n_0$ and $x\leq x_0$:
\begin{equation}
0\leq f_n(x)<\epsilon/3 \mbox{ and } 0\leq f(x)\leq\epsilon/6.
\end{equation}
This yields:
\begin{equation}\label{proofthpolya2}
|f_n(x)-f(x)|<\epsilon, \forall x\leq x_0, \; \forall n\geq n_0.
\end{equation}
Similarly, we can show that there exists $x_1>x_0$ and $n_1 \in \mathbb{N}$ such that:
\begin{equation}
|f_n(x)-f(x)|<\epsilon, \forall x\geq x_1, \; \forall n\geq n_1.\label{proofthpolya3}
\end{equation}
We can directly apply Theorem \ref{polyatheorem} to the sequence $(f_n)$ restricted to $[x_0,x_1]$, which gives the existence of $n_2 \in \mathbb{N}$, such that:
\begin{equation}
|f_n(x)-f(x)|<\epsilon, \forall x\in [x_0,x_1], \; \forall n\geq n_2.\label{proofthpolya4}
\end{equation}
By taking $n_3=\max \lbrace n_0, n_1, n_2 \rbrace$, and combining \eqref{proofthpolya2}, \eqref{proofthpolya3} \eqref{proofthpolya4} we get:
\begin{equation}
|f_n(x)-f(x)|<\epsilon, \forall x\in \mathbb{R}, \; \forall n\geq n_3,
\end{equation}
which concludes the proof.
\end{proof}

\end{appendices}
\section*{Acknowledgements}
I am very grateful to my thesis advisor, Julien Berestycki, for his help throughout the writing of this article. I would also thank Joon Kwon for very useful comments and the referees for giving me very precious advices.

 \bibliographystyle{plain}
 \bibliography{Biblio}
\end{document}